\input amstex
\documentstyle{amsppt}
\magnification=\magstep1

\pageheight{9.0truein}
\pagewidth{6.5truein}

\input xy
%%\input xypic

%%\xyoption{ps}
%%\UsePSspecials{Textures}

\xyoption{curve}\xyoption{arrow}\xyoption{matrix}

\def\edge{\ar@{-}}

\def\uloopr#1{\ar@'{@+{[0,0]+(-4,5)} @+{[0,0]+(0,10)} @+{[0,0]+(4,5)}}
  ^{#1}}
\def\dloopr#1{\ar@'{@+{[0,0]+(-4,-5)} @+{[0,0]+(0,-10)} @+{[0,0]+(4,-5)}}
  _{#1}}
\def\rloopd#1{\ar@'{@+{[0,0]+(5,4)} @+{[0,0]+(10,0)} @+{[0,0]+(5,-4)}}
  ^{#1}}
\def\lloopd#1{\ar@'{@+{[0,0]+(-5,4)} @+{[0,0]+(-10,0)} @+{[0,0]+(-5,-4)}}
  _{#1}}
\def\udotloopr#1{\ar@{{}.>} @'{@+{[0,0]+(-4,5)} @+{[0,0]+(0,10)}
  @+{[0,0]+(4,5)}} ^{#1}}

\def\urloop#1{\ar@'{@+{[0,0]+(1,9)} @+{[0,0]+(10,10)} @+{[0,0]+(9,1)}}
  ^{#1}}
\def\drloop#1{\ar@'{@+{[0,0]+(1,-9)} @+{[0,0]+(10,-10)} @+{[0,0]+(9,-1)}}
  _{#1}}
\def\ruloop#1{\ar@'{@+{[0,0]+(9,1)} @+{[0,0]+(10,10)} @+{[0,0]+(1,9)}}
  _{#1}}
\def\rdloop#1{\ar@'{@+{[0,0]+(9,-1)} @+{[0,0]+(10,-10)} @+{[0,0]+(1,-9)}}
  ^{#1}}

%%\long\def\ignore#1{}
\long\def\ignore#1{#1}

\def\seq{\mathrel{\widehat{=}}}
\def\la{{\Lambda}}
\def \len{\operatorname{length}} 
\def\detour{\mathrel{\wr\wr}}
\def\genus{\operatorname{genus}}

\def\dim{\operatorname{dim}}
\def\lamod{\Lambda\operatorname{-mod}}
\def\im{\operatorname{Im}}

\topmatter
\title The geometry of uniserial representations of finite dimensional
algebras I
\endtitle
\rightheadtext{Geometry of uniserial representations}
\author Birge Zimmermann Huisgen\endauthor
\address Department of Mathematics, University of California, Santa Barbara,
CA 93106\endaddress
\email birge\@math.ucsb.edu\endemail
\thanks This research was partially supported by a National Science
Foundation grant.\endthanks

\abstract It is shown that, given any finite dimensional, split basic algebra $\Lambda = K\Gamma/I$ (where $\Gamma$ is a quiver and $I$ an admissible ideal in the path algebra $K \Gamma$), there is a finite list of affine algebraic varieties, the points of which correspond in a natural fashion to the isomorphism types of uniserial left $\Lambda$-modules, and the geometry of which faithfully reflects the constraints met in constructing such modules. A constructive coordinatized access to these varieties is given, as well as to the accompanying natural surjections from the varieties onto families of uniserial modules with fixed composition series. The fibres of these maps are explored, one of the results being a simple algorithm to resolve the isomorphism problem for uniserial modules. Moreover, new invariants measuring the complexity of the uniserial representation theory are derived from the geometric viewpoint. Finally, it is proved that each affine algebraic variety arises as a variety of uniserial modules over a suitable finite dimensional algebra, in a setting where the points are in one-one correspondence with the isomorphism classes of uniserial modules. \endabstract

\endtopmatter
\document

\head 1. Introduction \endhead

This is the first part of a trilogy which is to lay the foundations for a
geometric approach to the uniserial representations of finite dimensional
algebras. The need for a solid understanding of the full class of uniserial
representations arose within a program aimed at approximating finitely
generated modules by modules of a simpler structure, the basic building blocks
of `helpful' approximations being uniserial modules. Of course, a uniserial
module viewed by itself does not hold much interest, as is the case for a point
on a curve when it is considered outside the context of the curve. In
classifying families of uniserial modules, the interest lies in the number of
and interplay among the parameters which offer themselves for the description
of these modules.

 For good focus, let us start by recording two problems which
have been propagated by Maurice Auslander since the 1970's and have now
appeared among the eleven open problems stated in [2, pp\. 411-412]:

(1) Give a method for deciding when two uniserial modules over an artin
algebra are isomorphic.

(2) Which artin algebras of infinite representation type have only a finite
number of pairwise nonisomorphic uniserial modules?

Specializing to finite dimensional algebras over algebraically
closed fields, we present
a solution to problem 1 and supply the tools to answer question 2 in a
subsequent article [7]. More generally, our results address split basic
finite dimensional algebras $\la$ over an arbitrary field $K$, i.e., algebras
of the form $K\Gamma/I$, where $\Gamma$ is a quiver and $I$ an admissible
ideal in the path algebra $K\Gamma$. Given such a `coordinatization' of $\la$,
we introduce irreducible affine algebraic varieties over $K$, the points of
which parametrize the isomorphism types of uniserial $\la$-modules,
and the geometry of which reflects the constraints met in constructing such
modules. While being {\it a priori} defined in terms of quiver and relations,
this list of varieties will turn out to be uniquely determined by the
isomorphism type of $\la$, up to order and birational equivalence [6].

There are several `classical' ways of viewing representations as points of
algebraic varieties; all of these varieties contain classes of uniserial
representations as open subvarieties. The most time-honored procedure is to
fix a $K$-basis for $\la$, say $\lambda_1,\dots,\lambda_m$, and to define the
variety $\operatorname{\Bbb Mod}_d(\la)$ of $d$-dimensional
$\la$-modules as a closed subvariety of $M_d(K)^m$ where the matrix in the
$i$-th slot represents module multiplication by $\lambda_i$. Obviously
these varieties are enormous -- as are those obtained by viewing
$d$-dimensional modules as points on Grassmannians; they include enough
redundancy to obscure potential connections between geometry and
representation theory. The sets of points corresponding to the isomorphism
classes of modules arise as orbits under the obvious
$GL_d(K)$-action (c.f. [8, \S\S12.16, 12.17]). It is not
difficult to see that the
$d$-dimensional uniserial modules with a fixed sequence of consecutive
composition factors form an open subvariety of
$\operatorname{\Bbb Mod}_d(\la)$, which is invariant under the canonical
$GL_d(K)$-action.  However, since the $GL_d(K)$-orbits can be explicitly
described only once the general isomorphism problem has been solved, these
varieties are not a useful resource for us. 

In contrast, the varieties which we will introduce provide a very close fit
to the families of uniserial modules they describe. In many cases,
there is a 1--1 correspondence between points and isomorphism types, and when
there is not, the `slack' occurring in the varieties is fairly harmless. 
In fact, there is a manageable way of calculating the fibres of the
correspondence in general. For more detail, we concentrate on the class of
uniserial representations of length $l+1$ having a fixed sequence
$\Bbb S= (S(1),\dots, S(l+1))$ of simple composition factors. It turns out
that there is a natural subdivision of this class, possibly with overlaps, so
that each of the segments is described by an {\it irreducible} affine
variety. A primary, somewhat rougher, subdivision is in terms of `masts', as
follows. Clearly, for each uniserial $\la$-module $U$ with sequence $\Bbb
S$ of consecutive composition factors, there exists at least one path $p$ of
length
$l$ in
$K\Gamma$ such that $pU\ne
0$; necessarily $p$ passes in order through the sequence $(e(1),\dots,e(l+1))$
of those vertices in $\Gamma$ which represent the simple modules $S(i)$. (Here
we consider $U$ as a $K\Gamma$-module in the obvious fashion.) Each such path
will be called a {\it mast} for $U$, and for each mast $p$ we will
construct an affine variety $V_p$ over $K$ -- not necessarily irreducible --
and a canonical surjection
$\Phi_p$ from $V_p$ onto the set of isomorphism types of uniserial left
$\la$-modules with mast $p$ (Section 3). If $\Gamma$ does not have any double
arrows, there is clearly at most one path $p$ passing through the sequence
$(e(1),\dots,e(l+1))$ of vertices, and the variety $V_p$ parametrizes the
full set of isomorphism classes of uniserial modules with composition sequence
$\Bbb S$. In case $\Gamma$ does contain double arrows, several varieties of
the form
$V_p$ are needed to account for the uniserial modules with composition
sequence $\Bbb S$ in general. In this case, the irreducible components of all
the pertinent $V_p$'s are combined into a variety $V_{\Bbb S}$, as discussed
below. Nonetheless, the `packaging' of irreducible components in terms of
masts remains the most accessible for both proofs and computations and hence
will play a crucial role in our work.

Loosely speaking, the points of
$V_p$ are -- as in the classical approach -- strings of coordinate
vectors determining module multiplication.  But the bases
for the (l+1)-dimensional
$K$-space underlying the representations in the image of $\Phi_p$ used here are
being shifted from one point of $V_p$ to the next, in a fashion that is tied
to the multiplication of the uniserial modules labeled by these points.  While
we lose the natural $GL_{l+1}$-action coming with
$\operatorname{\Bbb Mod}_{l+1}(\la)$ in the process, we gain, among other
things, a particularly transparent connection between the points of $V_p$ and
the graphs of the uniserial modules they represent, as well as a geometric
picture which cleanly shows how the relations in the ideal $I$ impinge on the
interplay of the parameters of the uniserials.  As another direct bonus of
restricting our attention to uniserial representations from the start, we can
make do with a small selection of paths from $K\Gamma$ as representative
`multipliers' and with bases for small $K$-subspaces of the modules considered
in order to pin down the effect of multiplication.  Unlike the classical
constructions, our sets of coordinate strings do not advertize themselves as
varieties at first sight.  We will actually start by describing a collection
of polynomials and subsequently prove that their vanishing locus
consists precisely of the mentioned families of parameters. The main price we
pay for the tight fit and manageability of the resulting varieties lies in the
fact that their independence of the chosen coordinate system for
$\la$, namely $\Gamma$ and $I$, is not at all obvious.  For instance, the
irreducible components of some
$V_p$ may shift to a different variety $V_q$ under a change of
coordinatization.  However, if we denote by
$V_{\Bbb S}$ the full collection of all the irreducible components of the
varieties
$V_q$, where $q$ runs through the paths of length $l$ that pass through the
sequence
$(e(1), \dots , e(l+1))$, the isomorphism type of $\la$ uniquely determines
$V_{\Bbb S}$ up to order and birational equivalence of the components.  In
case $\Gamma$ has no double arrows, the irreducible varieties in
$V_{\Bbb S}$ are even unique up to isomorphism; this is the situation in
which we find ourselves in case there are only finitely many ismorphism
types of uniserial $\la$-modules, for instance.  To understand the emerging
picture, we further require an in-depth study of intersections of the form
$\Phi_p(V_p)
\cap \Phi_q(V_q)$, where $p$ and $q$ are paths of length $l$ running through
the same vertices in the same order (Section 5). Since proving uniqueness and,
more generally, the functorial properties of our construction involves
considerable technical ballast, we address these issues separately in [6].

In case the path $p \in K\Gamma$ does not start with an oriented cycle, the
surjection 
$$\Phi_p : V_p \rightarrow \{\text{isomorphism types of uniserials
in\ } \lamod \text{with mast\ } p \}$$  
is even a bijection;  as a
consequence, the points of $V_p$ serve as complete families of isomorphism
invariants for the uniserial modules with mast $p$ in this situation.  On the
other hand,
$\Phi_p$ may fail to be bijective in the presence of certain types of
oriented cycles within $p$.  In Section 4, we characterize the fibres of the
maps $\Phi_p$, thus providing a general solution to the isomorphism problem. 
It turns out that there is a system of equations with the following property:
on insertion of arbitrary points $P$ and $Q$ of
$V_p$, it specializes to a {\it linear} system in the remaining variables, the
consistency of which is equivalent to `$\Phi_p(P) \cong
\Phi_p(Q)$'.  This
classification is quite gratifying since the varieties $V_p$, as well as the
pertinent systems of equations, are very accessible  --  they can be
readily computed on the basis of  quiver and relations for $\la$  --   and the
points
$P$ of
$V_p$ store information on the modules $\Phi_p(P)$ in a compact and easily
decodable form.  In [7], we will
further exploit this procedure from a theoretical point of view. In particular,
we will use it to tighten the connection between the fibres of
$V_p$ and certain oriented cycles of
$\Gamma$ and to solve the second of the problems 
stated at the outset.     

In Section 6, finally, we show that each affine variety over $K$ can be
realized as a variety $V_{\Bbb S}$
 for a suitable finite dimensional algebra $\la = K\Gamma/I$ and a suitable
sequence $\Bbb S$ of simple $\la$-modules.  We can even impose the
following additional requirements on the quiver $\Gamma$:  Namely, (a) that it
be without double arrows, which implies that $V_{\Bbb S}$ can be identified
with a single variety $V_p$, uniquely determined up to isomorphism by
$\la$, and (b) that it be acyclic, the latter condition entailfing bijectivity
of the corresponding natural map $\Phi_p$ as discussed above.  

The author would like to thank Axel Boldt for numerous helpful comments on a
first draft of this paper.

\head  2. Preliminaries \endhead 

	Throughout, $K$ will stand for an arbitrary field and $\la\cong K\Gamma/I$
will be a finite dimensional path algebra modulo relations over $K$; here
$\Gamma$ is a quiver and $I$ is an admissible ideal in the path algebra
$K\Gamma$. Our convention for the composition of paths $p,q\in K\Gamma$ is as
in [2], namely, $qp$ stands for `$q$ after $p$' whenever the concatenation is
defined. Moreover, $J$
will denote the Jacobson radical of $\la$. As is well known, the ideal
$I$ factored out of
$K\Gamma$ is, in general, not an isomorphism invariant of $\la$, but depends
on the choice of a complete set $e_1, \dots, e_n$ of primitive idempotents and
on the choice of `arrows' from $e_i$ to $e_j$, that is, of elements
from $e_jJe_i$ which give rise to a $K$-basis for $e_j
J e_i/ e_j J^2 e_i$. Any choice of an ideal $I$ in $K\Gamma$, together with
an isomorphism from $K\Gamma/I$ onto $\la$, will be called a {\it
coordinatization} of
$\la$.  For simplicity, we will assume that $\la$ is equal to $K\Gamma/I$,
unless otherwise specified, and identify the vertices of $\Gamma$ with the
corresponding primitive idempotents $e_1,\dots,e_n$ of $\la$. Our $\la$-modules will be {\it
left} modules throughout.   

\definition {Definition 1}  Given a uniserial $\la$-module  $U$ of length
$l+1$, any path
$p$ of length $l$ in  $K\Gamma$ with $pU \ne 0$ is called a {\it mast} of  $U$.
\enddefinition

\remark {Elementary Observations}  1. If the quiver $\Gamma$ has no double
arrows, then each uniserial $\Lambda$-module has a unique mast.  Conversely,
the uniqueness of masts in all uniserial modules implies absence of
double arrows. 

2. Not every path in $K\Gamma$ with nonzero image in $\la$ needs to occur as a
mast of a uniserial module. For example, if $\la= K\Gamma/I$, where $\Gamma$
is the quiver

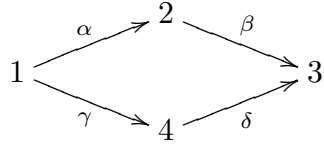
\ignore{
$$\xymatrixcolsep{3pc}\xymatrixrowsep{0.5pc}
\xy\xymatrix{
 & 2 \ar[dr]^\beta & \\
1 \ar[ur]^\alpha \ar[dr]_\gamma && 3 \\
 & 4 \ar[ur]_\delta & }\endxy$$
}

\noindent and $I= \langle \beta\alpha- \delta\gamma \rangle$, then neither
$\beta\alpha$ nor $\delta\gamma$ is a mast of a uniserial $\la$-module.

3. Of course, this concept of mast depends on a given coordinatization of
$\la$.  Our choice of coordinates will therefore impinge on the varieties of
uniserials with given mast, to be introduced in the next section.  As a
consequence, the effect of a coordinate change needs to be discussed (see
[6]).
\endremark
\medskip

 Let $M$ be a
$\la$-module. A {\it top element} of $M$ is an element $x\in M\setminus JM$
with $e_ix=x$ for some $i\in \{1,\dots,n\}$; in that case, $x$ will also be
called a {\it top element of type $e_i$}. Clearly, given a uniserial
$\la$-module $U$ of length $l+1$ with top element $x$, a path $p\in K\Gamma$
of length $l$ is a mast of $U$ if and only if $px\ne 0$. A useful tool in
visualizing and communicating uniserial modules will be their {\it
labeled and layered graphs}. Suppose that $p=\alpha_l \cdots\alpha_1$ where
each arrow $\alpha_i$ has starting point $e(i)$ and endpoint $e(i+1)$. The
labeled and layered graph of $U$ with respect to a top element $x$ and a mast
$p$ consists of (a) the mast $p$, drawn vertically with $e(1)$ at the top and
edges labeled by the arrows $\alpha_i$, together with (b) an edge labeled
$\omega$ from $e(i)$ to $e(j)$ whenever $i<j$ and $\omega$ is an arrow from
$e(i)$ to $e(j)$ such that $\la \omega\alpha_{i-1} \cdots\alpha_1x =J^{j-1}U$.

For example, let $\Gamma$ be the quiver

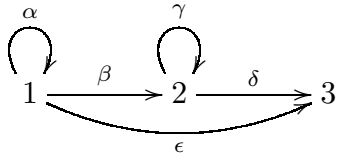
\ignore{
$$\xymatrixcolsep{3pc}
\xy\xymatrix{
1 \uloopr{\alpha} \ar[r]^\beta \ar@/_1pc/[rr]_\epsilon 
&2 \ar[r]^\delta \uloopr{\gamma} &3 
}\endxy$$
}

\noindent and $\la=K\Gamma/ \langle \alpha^2,\ \gamma^2 \rangle$. That the
(layered and labeled) graph of a uniserial module $U$ with top element $x$ be

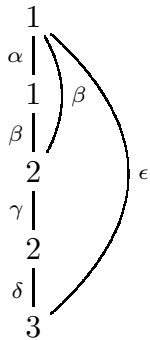
\ignore{
$$\xymatrixrowsep{1pc}
\xy\xymatrix{
1 \edge[d]_\alpha \edge @/^0.75pc/ [dd]^\beta
\edge @/^2.5pc/ [4,0]^\epsilon \\
1 \edge[d]_\beta \\
2 \edge[d]_\gamma \\
2 \edge[d]_\delta \\
3 }\endxy$$
}

\noindent means that $U$ has mast $p= \delta\gamma\beta\alpha$, that $\la
\beta x= J^2U$, and $\la \epsilon x= J^4U$; in other words, $\beta x$ is
congruent to a nonzero scalar multiple of $\beta\alpha x$ modulo $J^3U$, and,
$\epsilon x$ is a nonzero scalar multiple of $px$. Observe that these layered
and labeled graphs are, in general, not completely determined by the
isomorphism types of the modules they represent, but may depend on the choice
of top element. For instance, if $\la$ is as above and $U= \la e_1/ \bigl(
\la\gamma\beta\alpha +\la(\beta -\beta\alpha) +\la\epsilon \bigr)$, then the
graphs of $U$ relative to the top elements $x= \overline e_1$ and $y=
\overline e_1 -\alpha
\overline e_1$ are

\ignore{
$$\xymatrixrowsep{1pc}
\xy\xymatrix{
1 \edge[d]_\alpha \edge @/^1pc/ [dd]^\beta &&&&&& 1
\edge[d]_\alpha \\
1 \edge[d]_\beta &&& \txt{and} &&& 1 \edge[d]_\beta \\
2 &&&&&& 2 }
\endxy$$  
}

\noindent respectively. Conversely, layered and labeled graphs of uniserial
modules do not pin these modules down up to isomorphism. Indeed, if $\la=
K\Gamma$, where $\Gamma$ is the Kronecker quiver

\ignore{
$$\xymatrixcolsep{3pc}
\xy\xymatrix{
1 \ar[r]<0.5ex>^\alpha \ar[r]<-0.5ex>_\beta &2 }\endxy$$
}

\noindent then any uniserial module $U_k =\la e_1/ \la(\beta -k\alpha)$ for
$k\in K\setminus \{0\}$ has graph

\ignore{
$$\xymatrixrowsep{1.5pc}
\xy\xymatrix{
1 \edge[d]<-0.5ex>_\alpha \edge[d]<0.5ex>^\beta \\
2 }\endxy$$
}

\noindent while obviously $U_k \not\cong U_l$ for $k\ne l$.

Finally, we will speak of {\it subpaths} of a path $p\in K\Gamma$: A path
$q\in K\Gamma$ is a {\it right subpath} (respectively, {\it left subpath}) of
$p$ if there exists a path $r$ with $p=rq$ (respectively, $p=qr$). In
particular, if $p$ is a path from a vertex $e(1)$ to a vertex $e(2)$, then
$e(1)$ is the unique right subpath of length zero of $p$, and $e(2)$ is the
unique left subpath of length zero. A {\it proper} right subpath of $p$ is a
right subpath which is strictly shorter than $p$. It will be convenient to
communicate the statement `$q$ is a right subpath of $p$' in the form
`$p=\bullet q$' to arrive at compact formulas.

For the geometric terms we use, the reader should consult the introductory
texts by Hartshorne and Mumford ([5],[9]).

\head 3. Description of the varieties of uniserials with fixed mast \endhead

Throughout this section, let $p$ be a path of
length $l$ from a vertex $e(1)$ to a vertex $e(l+1)$.  We will use $p$ to
denote both the element in $K\Gamma$ and its residue class in $\Lambda$,
unless there is a danger of ambiguity.

Roughly speaking, the affine $K$-variety $V_p$ corresponding
to the uniserial $\Lambda$-modules with mast $p$ consists of points that are
families of coordinate vectors of the following type:  Given a uniserial left
$\la$-module
$U$ with mast $p$ and top element $x$, we equip $U$ with the $K$-basis
$p_ix$, for $p=\bullet p_i$, and string up the coordinate vectors of
the elements $qx$, where $q$ runs through the paths in $K\Gamma$ not vanishing
in $\la$.  Of
course, these coordinate strings pin down the corresponding uniserials.  The
first point to be addressed is the fact that these strings of coordinate
vectors actually form an affine algebraic variety $V_p$ which, in fact, is
readily accessible on the basis of a coordinatization of
$\la$.  A crude outline of the procedure for
assembling a finite set of polynomials which defines $V_p$ is as follows: 
Replace the scalars $k_{i,q}$ arising in the linear dependence relations $qx =
\sum_{p=\bullet p_i} k_{i,q} p_i x$ inside an arbitrary uniserial module with
mast
$p$ and top element
$x$ by generic indeterminates $X_{i,q}$, and expand a representative set of
relations from $I$ by means of substitutions $q= \sum_{p=\bullet p_i} X_{i,q}
p_i$ inside the polynomial ring
$K\Gamma[X_{i,q}]$.  A repetition of this substitution process will eventually
reduce the relations to equations of the form
$\sum_{p=\bullet p_i} \tau_i p_i = 0$, where the 
$\tau_i$ are polynomials in $K[X_{i,q}]$.  Reflecting the fact that
the elements $p_i x$ are $K$-linearly independent, our
interest will be in the simultaneous vanishing set of these polynomials
$\tau_i$. What makes the resulting affine variety fairly
manageable is the fact that we can significantly reduce the set of variables
without renouncing information.

Next we give a formal description of the variety $V_p$ and the canonical map
$\Phi_p$ from $V_p$ onto the set of isomorphism classes of uniserial
$\Lambda$-modules with mast $p$. Instead of considering all the
indeterminates $X_{i,q}$ mentioned above, we will restrict our attention to
those of the form $X_{i,\alpha u}$, where $(\alpha,u)$ is a `detour' as
defined below.

\definition {Definitions 2} 
 
{\bf 1.}  A {\it detour on the path} $p$ is a pair $(\alpha, u)$, where
$\alpha$ is an arrow and $u$ is a right subpath of $p$ (length 0 being
allowed) such that

(i) $\alpha u \ne 0$ in $K\Gamma$,

(ii) $\alpha u$ is not a right subpath of $p$ in $K\Gamma$, but 

(iii) there exists a right subpath $v$ of $p$ with $\len (v)\geq
\len (u) +1$ such that the endpoint of $v$ coincides with the endpoint of
$\alpha$.  

(See Figure (*) below.)  For ease of notation, we will often abbreviate the
statement `$(\alpha, u)$ is a detour on $p$' by `$(\alpha, u) \detour p$'.

{\bf 2.} Suppose that $p$ is a path of length $l$ that passes consecutively
through the vertices
$e(1),\dots,e(l+1)$. A {\it route on $p$} is any path in $K\Gamma$ which
starts in $e(1)$ and passes through a subsequence of the sequence
$(e(1),\dots,e(l+1))$. (We include $e(1)$ in the set of routes on $p$.)

\enddefinition

\remark{Remarks} (a) In particular, each right subpath of
$p$ is a route on $p$. More generally, any right subpath of a route is again a
route.

(b) In our work with routes, the following factorization property will be
crucial: A path $\bold r \in K\Gamma$ is a route on $p$ if and only
if it  can be written in the form  
$$\bold r = \bold r' \alpha _m u_m \cdots \alpha_1 u_1$$ 
for some $m \geq 0$
such that there exists a corresponding factorization 
$$p = p' \cdot w_m \cdots w_1$$ of
$p$ with the property that $(\alpha_i, u_i)$ is a detour on $w_i$ with
endpoint($\alpha_i$) = endpoint($w_i$) for each $i \leq m$, and $\bold
r'$ is a right subpath of $p'$.  (See Figure (**) below.)
Note that such a factorization of $p$ corresponding to the given
route $\bold r$ need not be unique.  

(c) The length of any route $\bold
r$ on
$p$ is bounded above by $\len(p)$.  

(d) Whenever $\bold r w$ is a route on $pw$, then $\bold r$ is a route on
$p$.\endremark

\ignore{
$$\xymatrixrowsep{1.5pc}
\xy\xymatrix{
{} \ar[dd]_u \ar[dd]<1.5ex>^u \ar@^{|-|}[4,0]<-3ex>_v
\ar@^{|-|}[6,0]<-6ex>_p &&&& {} \ar[dd]_{w_1} \ar[d]<1.5ex>^{u_1} 
\ar@^{|-|}[6,0]<-4ex>_p \\
 &&&& \ar @/^0.5pc/ [d]<1.4ex>^{\alpha_1} \\
{} \ar[dd] \ar @/^0.75pc/ [dd]<1.5ex>^\alpha  &&&& {} \ar[dd]_{w_2}
\ar[d]<1.5ex>^{u_2} \\
 &&&& \ar @/^0.5pc/ [d]<1.4ex>^{\alpha_2} \\
{} \ar[dd]   &&&& {} \ar[dd]_{p'} \ar[d]<1.5ex>^{\bold r'} \\
 &&&& \\
 &&&& \\
\txt{Figure (*)} &&&&\txt{Figure (**)}    }
\endxy$$
}

We start by describing the polynomial ring in which we will be working. Given
any detour $(\alpha,u)$ on $p$, let 
$$V(\alpha,u) =\{ v_i(\alpha, u)\mid i \in I(\alpha, u)\}$$
be the family of
right subpaths of $p$ in
$K\Gamma$ which are longer than $u$ and have the same endpoint as $\alpha$.
In most cases, it will be more convenient to refer to the index set
$I(\alpha,u)$ for $V(\alpha,u)$ than to the set $V(\alpha,u)$ itself.  The
polynomial ring of our choice will then be 
$$K\Gamma[X] =K\Gamma [X_i (\alpha, u)
\mid i \in I(\alpha, u), \ (\alpha, u)\detour p]$$
with coefficients in the path algebra $K\Gamma$ and independent variables
$X_i(\alpha,u)$. Next we introduce an equivalence relation on $K\Gamma[X]$ as
follows: Let
$\Cal L(p)$ be the left ideal of $K\Gamma[X]$ generated by all
the paths $q$ in
$K\Gamma$ which fail to be routes on $p$ together with all the differences
$$\alpha u - \sum_{i \in
I(\alpha, u)} X_i(\alpha, u) \cdot v_i(\alpha, u)$$
for detours $(\alpha,u)$ on $p$. Then,
clearly, the relation `$\sigma \seq \tau\ \Longleftrightarrow\ \sigma-\tau
\in \Cal L(p)$' for $\sigma,\tau \in K\Gamma[X]$ defines a congruence relation
relative to addition and left multiplication. The proof of the following easy
observation is left to the reader.

\proclaim{Observation 3} Each element of the path algebra $K\Gamma$ is
$\seq$-congruent to a unique element of the form $\sum_{p=\bullet p'}
\tau_{p'} (X) p'$, where the $\tau_{p'}(X)$ are polynomials in
$$K[X] = K[X_i (\alpha, u) \mid i \in I(\alpha, u), \ (\alpha, u)\detour
p].\qed$$
\endproclaim

To obtain these polynomials $\tau_{p'}(X)$ for a given element $z\in K\Gamma$
algorithmically, consider the following {\it substitution equations for} $p$:
First, $q \seq 0$ for any path $q$ in $K\Gamma$ which fails to
be a route on
$p$, and second,
$$\alpha u \seq \sum_{i \in
I(\alpha, u)} X_i(\alpha, u) \cdot v_i(\alpha, u)$$
for all detours $(\alpha, u)$ on $p$. We use the phrase `inserting the
substitution equations from the right' for the following steps:
$$\alignat2 q''q' &\seq 0 &&\qquad\qquad \text{if $q'$ is a non-route on
$p$};\\ 
q''q' &\seq \sum_{i \in I(\alpha, u)} X_i(\alpha, u) \cdot q''v_i(\alpha, u)
&&\qquad\qquad \text{if\ } q'=\alpha u \text{\ with\ } (\alpha,u)\detour p.
\endalignat$$  
Note that if $q'$ is not a route on $p$, then neither is
$q''q'$. Inserting the substitution equations from the right into the paths
occurring in the element
$z\in K\Gamma$ and repeating this procedure clearly leads to the equivalence
$z\seq
\sum
\tau_{p'}(X)p'$ with
$\tau_{p'}(X)\in K[X]$ after at most $d$ steps, where $d$ is an upper bound on
the lengths of the paths involved in $z$.

Now let $L$ be the Loewy
length of $\la$ -- i.e., $L$ is minimal with respect to $J^L = 0$ -- and denote
by $I^{(L)}$ the $K$-subspace of $I$ consisting of all elements which
can be written as $K$-linear combinations of paths of lengths at most $L$. 
Moreover, choose a finite $K$-generating set $t_1, \dots, t_s$ for the space
$I^{(L)}$.  By the above remarks, there are unique polynomials
$\tau_{i,p'}(X) \in K[X]$ with the property that $t_i \seq \sum_{p=\bullet p'}
\tau_{i,p'}(X) p'$ for $1 \leq i \leq s$. We now give a definition of the
variety $V_p$ depending on this choice of relations $t_1,\dots,t_s$; this is
the most convenient description for purposes of computation. However, as we
will immediately observe, this definition of $V_p$ is independent of
the choice of the $t_i$.

\definition {Definition 4}  Let $V_p = V( \tau_{i,p'}(X)\mid 1\le i\le s
\text{\ and\ } p=\bullet p')$ be the simultaneous vanishing locus of the
polynomials
$\tau_{i,p'}(X)$ in affine $N$-space $\Bbb A^N =\Bbb A^N(K)$, where $N = \sum
_{(\alpha,u) \detour p} |I(\alpha,u)|$.
\enddefinition

\remark{ Remarks}  
{\bf 1.} The affine algebraic set $V_p$ is independent of
the choice of a $K$-generating set for the $K$-space $I^{(L)}$; in
fact, $V_p$ is the vanishing locus of {\it all} polynomials $\tau_{p'}(X)$
arising in congruences $z \seq \sum_{p=\bullet p'}
\tau_{p'}(X) p'$ for some element $z$ from the ideal $I$.  The former
assertion is an immediate consequence of our construction, while the latter
can easily be checked as follows:  If the length of $p$ exceeds the Loewy
length
$L$ of $\la$, then $V_p$ is the empty set, because the right subpath $p_L$
of $p$ of length $L$ belongs to
$I^{(L)}$, and hence the congruence $p_L\seq 1 \cdot p_L$ places the constant
$1$ into the ideal of $V_p$.  If, on the other hand, $p$ has length $\leq L$,
all paths in $K\Gamma$ of length greater than $L$ are non-routes on $p$ and are
hence  reduced to zero under our equivalence relation.  

{\bf 2.} Our next remark along this line often saves a considerable amount of
computational effort. Namely, observe that, whenever $r_1,\dots,r_m\in I$
generate $I$ as a {\it left} ideal of $K\Gamma$, then $V_p= V(\rho_{i,q}(X)
\mid 1\le i\le m,\ p=\bullet q)$, where $r_i \seq \sum_{p=\bullet q}
\rho_{i,q}(X)q$ with $\rho_{i,q}(X)\in K[X]$.

{\bf 3.} A priori, the family of algebraic varieties $V_p$, where $p$ runs
through the set of paths  in $K\Gamma$ which do not vanish in $\la$, clearly
depends on the chosen coordinatization of the split algebra $\la$, that is, on
a fixed  set $ e_1, \dots, e_n$ of primitive idempotents and $K$-bases for the
spaces $e_i J e_j / e_i J^2 e_j$.  In fact, both the labels of the varieties
considered and their realization in affine space depend on the choice of
coordinates. This leads us to the following uniqueness problem: 
If
$K\Gamma / I \cong K\Gamma / I'$, are the varieties $V_p$ and $V'_p$
formed relative to the two ideals $I$ and $I'$ isomorphic? While in
general there is uniqueness only up to birational equivalence (see Section 5),
the answer to the  isomorphism question is positive for large
classes of algebras (see Section 6 and [6]).
\endremark      
 	
It is easy to compute
the varieties $V_p$ from a given coordinatization, i.e., from a 
presentation of
$\la$ in terms of quiver and relations.  These varieties will permit
us to classify the uniserial modules in terms of the correspondence described
in the main theorem of this section (Theorem A).  Before we state this theorem,
we illustrate the construction of the varieties $V_p$.

\example{Example 5} Let $\la= K\Gamma/I$, where $\Gamma$ is the quiver

\ignore{
$$\xymatrixcolsep{3pc}\xy\xymatrix{
1 \ar[r]<0.5ex>^{\alpha_1} \ar[r]<-0.5ex>_{\beta_1} 
&2 \ar[r]<0.5ex>^{\alpha_2} \ar[r]<-0.5ex>_{\beta_2} \uloopr{\gamma}
&3 \ar[r]<0.5ex>^{\alpha_3} \ar[r]<-0.5ex>_{\beta_3}
&4 \ar[r]<0.5ex>^{\alpha_4} \ar[r]<-0.5ex>_{\beta_4} &5
}\endxy$$
}

\noindent and $I$ is the ideal of $K\Gamma$ generated by the following
relations: $\gamma^3$, $\beta_2\alpha_1$, $\beta_2\gamma\alpha_1$,
$\alpha_4\alpha_3\beta_2\beta_1$, $\alpha_4\alpha_3\beta_2\gamma^2\alpha_1
-\alpha_4\alpha_3\alpha_2\gamma^2\alpha_1$,
$\beta_4\beta_3\alpha_2\gamma\alpha_1 -\alpha_4\alpha_3\alpha_2\beta_1$,
$\alpha_4\alpha_3\alpha_2\gamma^2\beta_1 -\alpha_4\alpha_3\alpha_2\alpha_1$,
$\alpha_4\alpha_3\alpha_2\gamma\alpha_1
-\alpha_4\alpha_3\beta_2\gamma\beta_1$, $\alpha_4\alpha_3\alpha_2\gamma\alpha_1
-\alpha_4\beta_3\alpha_2\gamma^2\alpha_1$,
$\alpha_4\alpha_3\alpha_2\gamma\alpha_1
-\beta_4\alpha_3\alpha_2\gamma^2\alpha_1$. Moreover, consider the path $p=
\alpha_4\alpha_3\alpha_2\gamma^2\alpha_1$. To compute the variety $V_p$,
observe that the detours on $p$ and the corresponding substitution equations
are as follows:
$$\alignat2 (\beta_1,e_1) &\qquad\qquad& \beta_1 &\seq X_1\alpha_1
+X_2\gamma\alpha_1 +X_3\gamma^2\alpha_1\\
(\alpha_2,\alpha_1) &\qquad\qquad& \alpha_2\alpha_1 &\seq
X_4\alpha_2\gamma^2\alpha_1\\
(\alpha_2,\gamma\alpha_1) &\qquad\qquad& \alpha_2\gamma\alpha_1 &\seq
X_5\alpha_2\gamma^2\alpha_1\\
(\beta_2,\alpha_1) &\qquad\qquad& \beta_2\alpha_1 &\seq
X_6\alpha_2\gamma^2\alpha_1\\
(\beta_2,\gamma\alpha_1) &\qquad\qquad& \beta_2\gamma\alpha_1 &\seq
X_7\alpha_2\gamma^2\alpha_1\\
(\beta_2,\gamma^2\alpha_1) &\qquad\qquad& \beta_2\gamma^2\alpha_1 &\seq
X_8\alpha_2\gamma^2\alpha_1\\
(\beta_3,\alpha_2\gamma^2\alpha_1) &\qquad\qquad& \beta_3\alpha_2\gamma^2\alpha_1
&\seq X_9\alpha_3\alpha_2\gamma^2\alpha_1\\
(\beta_4,\alpha_3\alpha_2\gamma^2\alpha_1) &\qquad\qquad&
\beta_4\alpha_3\alpha_2\gamma^2\alpha_1 &\seq X_{10}p\endalignat$$
Next observe that the relations listed above, together with the paths
$\gamma^3\alpha_1$ and $\gamma^3\beta_1$, generate $I$ as a left ideal of
$K\Gamma$; consequently, Remark 2 following Definition 4 tells us that we
need only consider these elements of $I$ in determining a generating set of
polynomials for the ideal of $V_p$. Since these last two paths, as well as
$\gamma^3$, are non-routes on $p$, they are $\seq$-equivalent to 0 and hence
do not lead to conditions on the indeterminates $X_i$. We now insert the
substitution equations into the remaining relations. 

First $\beta_2\alpha_1 \seq X_6\alpha_2\gamma^2\alpha_1$,
which, in view of $\beta_2\alpha_1\in I$, gives us the equation $X_6=0$ for
$V_p$. Analogously, the combination of
$\beta_2\gamma\alpha_1
\seq X_7\alpha_2\gamma^2\alpha_1$ and $\beta_2\gamma\alpha_1 \in I$ implies
$X_7=0$ for the points of
$V_p$.  Moreover, $\alpha_4\alpha_3\beta_2\beta_1 \seq
\alpha_4\alpha_3\beta_2(X_1\alpha_1 +
X_2\gamma\alpha_1 +X_3\gamma^2\alpha_1) \seq X_1X_6p +X_2X_7p +X_3X_8p$
yields $X_1X_6 +X_2X_7 +X_3X_8=0$ on $V_p$. In view of $X_6= X_7 =0$, we
obtain $X_3X_8=0$.  We further compute $\alpha_4\alpha_3\beta_2\gamma^2\alpha_1
-p \seq (X_8 -1)p$ to conclude that $X_8-1 =0$, and consequently also $X_3=0$,
on $V_p$.

Next, $\beta_4\beta_3\alpha_2\gamma\alpha_1 -\alpha_4\alpha_3\alpha_2\beta_1
\seq X_5\beta_4\beta_3\alpha_2\gamma^2\alpha_1
-X_1\alpha_4\alpha_3\alpha_2\alpha_1
-\allowmathbreak X_2\alpha_4\alpha_3\alpha_2\gamma\alpha_1 -\allowmathbreak
X_3\alpha_4\alpha_3\alpha_2\gamma^2\alpha_1\allowmathbreak
\seq X_5X_9X_{10}p - X_1X_4p -X_2X_5p -X_3p$, which gives us $X_5X_9X_{10}
-X_1X_4 -X_2X_5 -X_3 =0$, and
$X_5X_9X_{10} -X_1X_4 -X_2X_5 =0$ in view of the preceding equations.

Analogously, the last four relations yield $X_1-X_4=0$, $X_5 -X_1X_7 -X_2X_8
=X_5 -X_2 =0$, $X_5-X_9 =0$, and $X_5 -X_{10} =0$, respectively. 

We deduce that
$$V_p = V( X_3, X_6, X_7, X_8-1, X_5X_9X_{10} -X_1X_4
-X_2X_5, X_1-X_4, X_2 -X_5, X_5-X_9, X_5-X_{10}).$$
Clearly, $V_p$ is isomorphic to the variety  $V(Y^2-X^3-X^2)$ in affine
2-space. Over the real numbers, we thus obtain the standard $\alpha$-curve:

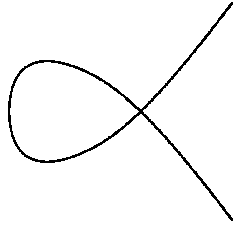
\ignore{
$$\xy (25,12);(25,-12)
\curve{(16,0)&(8,-6)&(0,-6)&(0,6)&(8,6)&(16,0)}
\endxy$$
}

\noindent In this example, the correspondence between the points of $V_p$ and
the isomorphism types of uniserial $\la$-modules with mast $p$, as described in
Theorem A below, is bijective, even though our path $p$ does
include an oriented cycle.\qed\endexample

The following theorem contains basic information about the `back and
forth' between points on the variety $V_p$ on one hand and uniserial
$\la$-modules with mast
$p$ on the other.

\proclaim {Theorem A}  Suppose that $p$ is a path in $K\Gamma$ starting in the
vertex $e(1)$.
 
{\rm (I)}  There is a surjective map $\Phi_p$ from the variety
$V_p$ to the set of isomorphism types of uniserial left $\la$-modules with mast
$p$.  It assigns to each point $k = \bigl(k_i (\alpha, u) \bigr)_{i \in I
(\alpha, u),\ (\alpha, u)\detour p}$ in $V_p$ the isomorphism type of the
module
$\la e(1)/ U_k$, where 
$$U_k = \Biggl(\sum_{(\alpha, u) \detour p} \la \biggl(\alpha u -
\sum_{i
\in I(\alpha, u)} k_i(\alpha, u) v_i (\alpha, u) \biggr)\Biggr) +
\Biggl(\sum_{q\
\text{not a route on }\ p} \la qe(1) \Biggr).$$  
Alternately, the uniserial module $\la e(1)/U_k$ representing $\Phi_p(k)$ can
be described as the unique uniserial factor module with mast $p$ of the module 
$$\la e(1) / \Biggl(\sum_{(\alpha,
u)\detour p} \la \biggl(\alpha u - \sum_{i \in I(\alpha, u)}
k_i(\alpha, u) v_i (\alpha, u) \biggr) \Biggr).$$ 

{\rm (II)}  The variety $V_p$ is nonempty if and only if there
exists a uniserial left $\la$-module with mast $p$.

{\rm (III)}  Provided that $p$ does not have a proper right subpath which is an
oriented cycle of positive length, the map
$\Phi_p$ is bijective. More strongly: If $p$ is not of the form $p'c$, where
$c$ is an oriented cycle and both $p'$ and $c$ are paths of positive length,
then $\Phi_p$ is bijective. 
\endproclaim

\demo {Proof of Theorem A} The following maps will repeatedly prove useful.  
For any family of scalars $k = \bigl(k_i (\alpha, u) \bigr)_{i \in I
(\alpha, u),\ (\alpha, u)\detour p}$ in $\Bbb A^N$, we will
consider a $K$-linear transformation $F_k: K\Gamma \rightarrow
K\Gamma$ depending on $k$:  Namely, $F_k$ sends all the paths which fail
to be routes on $p$ to zero and acts on routes $q$ as follows: If $q$
is a right subpath of $p$, then $F_k(q) = q$, and if $q = q'\alpha
u$, where $(\alpha, u)$ is a detour on $p$, then $F_k(q) = \sum_{i \in
I(\alpha,u)} k_i(\alpha, u) q' v_i(\alpha, u)$. Then, given any element
$z \in K\Gamma$ which is a $K$-linear combination of paths of length at
most $m$, the image of $z$ under $(F_k)^m$ equals $\sum_{p=\bullet p'}
\tau_{p'} (k) p'$, where $\sum_{p=\bullet p'} \tau_{p'} (X) p'$ is the
unique element in $\sum_{p=\bullet p'}K[X] p'$ equivalent
to $z$ under `$\seq$'; here we are working inside the polynomial ring
$K\Gamma [X]$, where $X$ stands for the family of variables $X_i(\alpha,u)$. 
In
particular, if $z$ belongs to $I$, the polynomials
$\tau_{p'} (X)$  fall into the ideal of $V_p$ (by the first of the above
remarks) and consequently vanish at $k$ whenever $k\in V_p$; in other words,
$(F_k)^m(z) =0$ in that case.

To any family $k = \bigl( k_i(\alpha, u) \bigr) \in \Bbb A^N$  as
above, we moreover assign the left ideal $W_k$ of the path
algebra
$K\Gamma$, generated by the paths $q$ which are non-routes on $p$ and by the
elements
$\alpha u - \sum_{i \in I(\alpha,
 u)} k_i(\alpha, u) v_i(\alpha, u)$, where $(\alpha, u)$ runs through
the detours on $p$.  Observe that, for any element $z \in K\Gamma$ and any
choice of $k$, the difference $z - F_k(z)$ belongs to $W_k$;  in
particular, the subspace $W_k$ of $K\Gamma$ is invariant under $F_k$.     

Let $\Phi = \Phi_p$ be defined as in the statement of the theorem and, for
$k\in V_p$, identify $\Phi(k)$ with the module $\la e(1)/U_k$.

{\it Claim 1.} For each $k \in V_p$, the module $\Phi (k)$ is
uniserial with mast $p$.

Indeed, write $p = \beta_l \cdots \beta_1$, where the $\beta_i$ are arrows,
and set $x = e(1) + U_k \in \Phi(k)$.  It clearly suffices to show that
$\Phi(k)$ has $K$-basis $\beta_i \cdots \beta_1 x$, $ 0\leq i \leq l$. 
Indeed, once this is established, we have $J^l \Phi(k) \ne 0$, and, for reasons
of dimension, the radical powers $J^i \Phi (k)$, $ 0 \leq i \leq l$,
constitute a composition series for $\Phi(k)$.

To see  that the elements $\beta_i \cdots \beta_1 x$, $ 0\leq i \leq l$, form a
$K$-generating set for $\Phi(k)$, it suffices to check that, given any
route $\bold r$ on $p$, the element $\bold r x\in \Phi(k)$ is a
$K$-linear combination of the $\beta_i \cdots \beta_1 x$.  We prove this by
induction on
$d = l -
\len (u_1)$, where $u_1$ is the right subpath of $\bold r$ of maximal length
occurring also as a right subpath of $p$.  If $d = 0$, then $\bold r =
p$, and there is nothing to prove. So suppose $d \geq 1$, and let $\bold r =
\bold r'\alpha_s u_s \cdots \alpha_1 u_1$.  If $s=0$, then $\bold r$ is a
(proper) right subpath of $p$, and we are again done. If, on the other
hand, $s\ge1$, the definition of $\Phi(k)$ yields  
$$\bold r x = \sum_{i \in I(\alpha_1, u_1)} k_i(\alpha_1, u_1)\bold r'
\alpha_s u_s \cdots \alpha_2 u_2 v_i(\alpha_1, u_1) x.$$  
If the path
$\bold r'\alpha_s u_s \cdots \alpha_2 u_2 v_i(\alpha_1, u_1)$ is not a route on
$p$, then $\bold r'\alpha_s u_s \cdots \alpha_2 u_2 v_i(\alpha_1, u_1) x = 0$
by construction.  The induction hypothesis applies to each of the remaining
terms
$\bold r'\alpha_s u_s
\cdots \alpha_2 u_2 v_i(\alpha_1, u_1) x$, since the lengths of the right
subpaths $v_i(\alpha_1, u_1)$ of $p$ exceed the length of $u_1$.

To prove that the elements $\beta_i \cdots \beta_1 x$ for $ 0\leq i \leq l$ are
linearly independent, it suffices to verify that $px = \beta_l \cdots \beta_1
x \ne 0$.  Assume the contrary, and consider the following presentation of
$\Phi (k)$ as a $K\Gamma$-module:  $\Phi (k) = K\Gamma e(1) /( Ie(1) +
W_ke(1))$, where
$W_k$ is as defined above.
So $px = 0$ is equivalent to the existence of an element $a \in I e(1)$ such
that $p + a \in W_k$. To reach a contradiction, we will first derive that
$p\in W_k$. Suppose that
$a$ is a
$K$-linear combination of paths of lengths at most $m$.  In view of the fact
that $a \in I$, the first paragraph of the proof shows that $F_k ^m(a) = 0$. 
Moreover, due to the
$F_k$-invariance of $W_k$, the fact that $p + a \in W_k$ yields 
$p = F_k^m(p + a) \in W_k$ as required. 
But the containment $p\in W_k$ is in turn impossible, as we now prove.

Assume that
$p$ is a
$K$-linear combination of terms  
$q'\bigl( \alpha u - \sum_{i \in I(\alpha, u)}k_i(\alpha, u)
v_i(\alpha, u) \bigr)$
and terms $q''q$, where $q$ is not a route on $p$ and $q'$,
$q''$ are paths in $K\Gamma$.  Since none of the paths of the form $\alpha u$
going with a detour $(\alpha, u)$ on $p$ is a right subpath of $p$, the above
equality forces $p$ to be equal to one of the terms $q' v_i(\alpha, u)$, say
$p= q'_0 v_i(\alpha_0,u_0)$. This clearly implies that
$q'_0\alpha_0u_0$ is a route on $p$. Observe that, more generally,
$q'\alpha u$ is a route on $p$, whenever $(\alpha,u)$ is a detour with the
property that
$q'v_j(\alpha,u)$ is a route on $p$ for some $j\in I(\alpha,u)$. Let $u_1$ be
a right subpath of $p$ that has minimal length with respect to the following
properties:

(a) for some arrow $\alpha_1$, the pair $(\alpha_1,u_1)$ is a detour on $p$,
and, for a suitable path $q_1'$, the element $q'_1 \bigl(\alpha_1u_1 - \sum
k_i(\alpha_1,u_1) v_i(\alpha_1,u_1) \bigr)$ occurs nontrivially in the above
representation of $p$ as an element of $W_k$;

(b) $q'_1\alpha_1u_1$ is a route on $p$.

Since $\alpha_1u_1$ is not a right subpath of $p$ and, a fortiori,
$q'_1\alpha_1u_1 \ne p$, the term $q'_1\alpha_1u_1$ must cancel out of our
representation of $p$ as an element of $W_k$. Clearly $q'_1\alpha_1u_1 \ne
q''q$ whenever $q$ fails to be a route on $p$, and $q'_1\alpha_1u_1 \ne
q'_2\alpha_2u_2$ whenever $(\alpha_2,u_2)$ is a detour different from
$(\alpha_1,u_1)$. Thus $q'_1\alpha_1u_1 =q'_2v_j(\alpha_2,u_2)$ for some $j$
and some detour $(\alpha_2,u_2)$ of $p$ such that $q'_2 \bigl(\alpha_2u_2 -
\sum_{i\in I(\alpha_2,u_2)} k_i(\alpha_2,u_2) v_i(\alpha_2,u_2) \bigr)$ occurs
non-trivially in our representation of $p$. On one hand, this forces
$q'_2v_j(\alpha_2,u_2)$, and consequently also $q'_2\alpha_2u_2$, to be a
route on $p$, whence $\len(u_2) \ge \len(u_1)$ by our choice of $u_1$. On the
other hand, the equality `$q'_1\alpha_1u_1 =q'_2v_j(\alpha_2,u_2)$' is only
possible when $\len v_j(\alpha_2,u_2)\le \len(u_1)$, for $u_1$ is the
longest right subpath of $p$ which is also a right subpath of
$q'_1\alpha_1u_1$. But this, in turn, entails $\len(u_2)< \len(u_1)$. We have
reached a contradiction, which proves that $p\notin W_k$, and which
thus completes the proof of Claim 1.

{\it Claim 2.}  Whenever $U$ in $\la$-mod is a uniserial module with
mast $p$, there exists a point $k \in V_p$ such that $U \cong
\Phi(k)$.

To find a suitable point $k \in V_p$, let $x = e(1)x$ be a top
element of $U$; then the map $f: \la e(1) \rightarrow U = \la x$
which sends $e(1)$ to $x$ is a projective cover of
$U$.  Again write $p = \beta_l \cdots \beta_1$, where the $\beta_i$
are arrows. Then, clearly, $J^iU = \la \beta_i \cdots \beta_1 x$ for $i
\le l$, and the elements $\beta_i \cdots \beta_1 x$, $0 \le i \le l$, form a
$K$-basis for
$U$. 

If $(\alpha, u)$ is a detour on $p$, then $\la \alpha u x = J^m x$ for some
$m \ge \len( u) + 1$.  Consequently, letting $e$ be the primitive
idempotent with $e\alpha = \alpha$, we obtain that $\alpha u x$ is
a $K$-linear combination of those elements $\beta_i \cdots \beta_1
x$ of $U$ for which $\beta_i \cdots \beta_1$ is a path longer than $u$
ending in $e$.  But these latter paths are precisely the ones
that belong to the family $\bigl( v_i(\alpha, u) \bigr)_{i \in
I(\alpha, u)}$, which shows that $\alpha u x = \sum_{i \in
I(\alpha, u)} k_i(\alpha, u) v_i(\alpha, u)x$ for suitable scalars
$k_i (\alpha, u) \in K$.  (Note on the side:  Once we have fixed a top
element $x$ of U, the scalars $k_i (\alpha, u)$ are actually uniquely
determined by $U$, since the elements $v_i(\alpha, u) x $, $ i \in I(\alpha,
u)$, are $K$-linearly independent.)  In other words, via the map $f$, the
uniserial $U$ is an epimorphic image of the module  $M := \la e(1)/
\sum_{(\alpha, u) \detour p} \la \bigl(\alpha u -  \sum_{i \in I(\alpha, u)}
k_i(\alpha, u) v_i (\alpha, u) \bigr)$.  Set $k = \bigl( k_i (\alpha, u)
\bigr)_{ i \in I(\alpha, u),\ (\alpha, u)\detour p}$.  In order to
prove Claim 2, it is thus enough to show that (a) $qx = 0$, whenever $q \in
K\Gamma$ is a path which fails to be a route on $p$, and (b) $k \in V_p$.  
For reasons of dimension, we will then obtain that the epimorphism
$\Phi(k) \rightarrow U$ induced by $f$ is an isomorphism. In particular,
statement (a) will guarantee that, up to isomorphism, $U$ is the only factor
module of $M$ which is uniserial with mast $p$.

For (a), suppose that the path $q$ is a non-route on $p$.  It is clearly
harmless to assume that $q = q e(1)$.  Write $q = q'u$, where  $u$ is the
longest right subpath of $q$ (possibly of length zero) which {\it is} a
route on $p$.  Moreover, let $p'$ be the shortest right subpath of $p$ such
that $u$ is a route on $p'$, say $p = p'' p'$. Then $\len q' \ge 1$, that is,
$q' =  q''\gamma$ for some arrow $\gamma$, and the fact that $\gamma u$ is
not a route on $p$ is equivalent to the non-existence of a right subpath $v$
of $p$ which is strictly longer than  $p'$ and ends in the same vertex as
$\gamma$.  If $ux \ne 0$, then the minimal choice of $p'$ implies that $\la
ux  =\la w'x$ where $w'$ is a right subpath of $p$ with $\len(w') \ge
\len(p')$. If $p=w''w'$, then $w''$ is a left subpath of
$p''$ such that
$\la ux =\la w'x$ is a uniserial module with mast $w''$. Consequently, the
fact that
$p''$, and hence also $w''$, is devoid of right subpaths of positive length
ending in the same vertex as $\gamma$ forces $\gamma ux$ to be zero. Thus $q x
=0$ as required.

To verify that $k$ is a point in $V_p$, let $z$ be any element in $I$; say $z$
is a linear combination of paths of lengths bounded above by some integer
$m$.  Moreover, let $\sum_{p=\bullet p'}
\tau_{p'}(X) p'$ be the unique element in $\sum_{p=\bullet p'}K[X] p'$ which
is equivalent to $z$ under $\seq$.  We want to show that all of the
polynomials $\tau_{p'} (X)$ vanish at $k$.  For that purpose, view
$U$ again as a $K\Gamma$-module under the action induced by that of
$\la$, and note that $zU = 0$, as well as $W_k x = 0$.  Since $F_k(y) - y \in
W_k$ for all $y \in K\Gamma$, we infer that $F_k^r (z) x = 0$ for all $r \ge
1$.  But as we pointed out in the first paragraph, $F_k ^m (z) = 
\sum_{p=\bullet p'}\tau_{p'}(k) p'$, and so, in
particular, $\sum_{p=\bullet p'}\tau_{p'}(k) p' x =
0$.  Now use the fact that the elements $p'x$, where $p'$ runs through the
right subpaths of $p$, are $K$-linearly independent, to conclude
$\tau_{p'}(k) = 0$ for all $p'$ as required.  This completes the proof of
part (I).  Part (II) is an obvious consequence of Part (I).

(III)  In  proving the surjectivity of $\Phi$, we saw that, given a
uniserial module $U$ with mast $p$ and a fixed top element $x$, there is a
unique point
$k = \bigl( k_i(\alpha, u) \bigr) \in \Phi ^{-1}(U)$ such that $\alpha u x =
\sum_{i \in I(\alpha, u)} k_i(\alpha, u) v_i(\alpha, u) x$ for each detour
$(\alpha, u)$ on $p$. In other words, each subfamily
$\bigl(k_i(\alpha, u)\bigr)_{i \in I(\alpha, u)}$ of $k$ is the 
coordinate vector of $\alpha u x$ relative to the $K$-basis $v_i(\alpha, u)
x$, $i\in I(\alpha,u)$, of $eJ^{\len(u)+1}U$; here $e$ is the primitive
idempotent in which
$\alpha$ terminates.  Since this coordinate vector will not change if $x$ is
replaced by $ax$ for some nonzero scalar $a$, we see that, if $U$ has a unique
top element, up to scalar factors  --  equivalently, if $p$ is not of the form
$p = p' c$ for a cycle $c$ of positive length from $e(1)$ to $e(1)$  --  the
set
$\Phi^{-1}(U)$ is a singleton, i.e., $\Phi$ is injective.  This leaves
us to deal with the case, where $p$ is a cycle from $e(1)$ to $e(1)$
such that no proper right subpath of $p$ is a cycle $e(1)\rightarrow e(1)$.  In
that case, an arbitrary top element
$y$ of
$U$ is of the form
$y = ax + bpx$ for $a,b \in K$ and $a$ nonzero, and since $\alpha u y = \alpha
u ax$ and
$v_i(\alpha, u) y = v_i(\alpha, u) ax$ for each detour $(\alpha, u)$ on $p$,
the coordinate vector of $\alpha u y$ with respect to the new basis is the
same as that of $\alpha ux$ with respect to the old basis.  Again we deduce
that
$\Phi$ is injective.
\qed\enddemo

Note that, in case $K$ is algebraically closed, Theorem A(II) provides us
with an algorithmic procedure to decide whether a given path $p :
e\rightarrow e'$ in $K\Gamma$ occurs as mast of a uniserial left
$\Lambda$-module. Indeed, when combined with Hilbert's Nullstellensatz,
Theorem A yields the following: There is a uniserial $\Lambda$-module with
mast $p$ if and only if, for some $K$-generating set $t_1,\dots,t_s$ for
$I^{(L)}$, the resulting polynomials $\tau_{ij}(X)$, as in the
definition of the variety $V_p$, generate a proper ideal of $K[X]$.
The ensuing question, whether or not 1 belongs to the ideal generated by the
$\tau_{ij}(X)$, $1\le i\le s$, $1\le j\le l$, is well-known to be decidable
by way of the Groebner method.

We will extract an observation from the preceding argument which will be
useful on several occasions. For a smooth formulation, we require the
following notational convention: Given a uniserial module $U$ with mast $p$
and a detour $(\alpha,u)$ on $p$ such that $\alpha$ ends in the vertex $e$,
we denote by $U(\alpha,u)$ the $K$-subspace $eJ^{\len(u)+1}U$ of $U$. Note
that, for any top element $x$ of $U$, the set $V(\alpha,u)x =\{v_i(\alpha,u)x
\mid i\in I(\alpha,u)\}$ forms a basis for this subspace.

\proclaim {Corollary to the proof of Theorem A}  Let $N$ be the disjoint
union of the index sets 
$I(\alpha, u)$, where $(\alpha, u)$ runs through the detours on $p$, and let
$k$ be a point in $\Bbb A^N(K)$, say $k = \bigl( k_i(\alpha, u) \bigr)_{i
\in I(\alpha, u),\ (\alpha, u)\detour p}$. Then $k$ belongs to $V_p$
precisely when there exists a uniserial module $U$ with mast $p$ and top
element $x$ such that, for each detour $(\alpha, u)$ on $p$, the projection
$\bigl( k_i(\alpha, u)
\bigr)_{i\in I(\alpha, u)}$ of $k$ onto $\Bbb A^{I(\alpha,
u)}$ is the coordinate vector of
the element $\alpha u x \in U(\alpha,u)$ with respect to the $K$-basis
$V(\alpha,u)$. 

In the positive case, any such uniserial module $U$ is isomorphic to
$\Phi_p(k)$, and the top element $x= e(1) +U_k$ of $\Phi_p(k)$ has the
property that $\alpha ux =\sum_{i\in I(\alpha,u)} k_i(\alpha,u)
v_i(\alpha,u)x$ for all $(\alpha,u)\detour p$; here $U_k$ is as in the
statement of Theorem A.
\qed\endproclaim

Let $L$ again denote the Loewy length of $\la$, and $I^{(L)}$ the $K$-subspace
of $I$ consisting of all elements of $I$ which can be written as $K$-linear
combinations of paths of lengths $\le L$. In the definition of $V_p$, we
picked a $K$-generating set for $I^{(L)}$ to arrive -- via the substitution
equations for $p$ -- at a set of polynomials that determines $V_p$. This set
may be vastly redundant, as we already pointed out after Definition 4. In
fact, it suffices to consider elements
$t_1,\dots,t_r\in I^{(L)}$ such that $I^{(L)} \subseteq \sum_{i=1}^r K\Gamma
t_i$; if again
$\tau_{ij}(X)\in K[X]$ are such that $t_i \seq \sum_{j=0}^l \tau_{ij}(X)p_j$
for $1\le i\le s$, then 
$$V_p= V(\tau_{ij}(X) \mid 1\le i\le r,\ 0\le j\le l).$$
On the other hand, it does not suffice to consider a set of relations which
generates
$I$ as an ideal, as the following example demonstrates.

\example{Example 6} Let $\la =K\Gamma/I$, where $\Gamma$ is the quiver

\ignore{
$$\xymatrixcolsep{3pc}
\xy\xymatrix{
1 \ar[r]^\alpha &2 \ar[r]<0.5ex>^\beta \ar[r]<-0.5ex>_\gamma 
&3 \ar[r]<0.5ex>^\delta \ar[r]<-0.5ex>_\epsilon &4
}\endxy$$
}

\noindent and $I\subseteq K\Gamma$ is the ideal generated by the
relations
$\delta\beta -\epsilon\gamma$, $\epsilon\beta$, and $\delta\gamma$. Moreover,
let $p= \delta\beta\alpha$. The detours on $p$ are $(\gamma,\alpha)$ and
$(\epsilon, \beta\alpha)$, and hence the substitution equations for $p$ are
$\gamma\alpha \seq X_1\beta\alpha$, $\epsilon\beta\alpha \seq X_2p$, as well
as
$q\seq 0$ whenever $q$ fails to be a route on $p$. Note that none of the
paths occurring in the above relations is a route on $p$, and so the
polynomials resulting from these relations via the substitution equations are
all zero. On the other hand, the relations $\delta\beta\alpha
-\epsilon\gamma\alpha$ and $\epsilon\beta\alpha$ yield $X_1X_2-1=0$ and
$X_2=0$, whence $V_p=\varnothing$. \qed\endexample

The following instance of failure of bijectivity of $\Phi_p$ is prototypical.

\example{Example 7} Let $\Lambda =K\Gamma/ \langle\alpha^2\rangle$, where
$\Gamma$ is the quiver

\ignore{
$$\xymatrixcolsep{3pc}
\xy\xymatrix{
1 \lloopd{\alpha} \ar[r]^\beta &2
}\endxy$$
}

\noindent If $p= \beta\alpha$, then
$V_p= \Bbb A^1$, while $\Phi_p(V_p)$ is a singleton; indeed, for any scalar
$k\in K$, the uniserial modules $\Lambda
e_1/\Lambda(\beta -k\beta\alpha)$ and $\Lambda e_1/\Lambda\beta$ are
isomorphic.

However, the sufficient condition for bijectivity of $\Phi_p$ given in
Theorem A(III) is not necessary.  Indeed, if $\Gamma$ is the quiver

QUIVER

and $\la = K\Gamma/\langle \alpha\gamma\alpha \rangle$, the path $p =
\beta\gamma\alpha$ again consists of a cycle followed by a nontrivial left
subpath.  Again we have $V_p= \Bbb A^1$, but this time the map
$\Phi_p$ is bijective.   
\qed\endexample

\definition{Definition 8}  Again, let $p$ be a path in $K\Gamma$ and $\la=
K\Gamma/I$. We will refer to $V_p$ as the {\it uniserial variety of
$\lamod$ at $p$}. The irreducible components of the uniserial variety
$V_p$ will be called the {\it uniserial components of $\lamod$ at $p$}.
Moreover, we will say that a uniserial component $W$ of $V_p$  {\it
intersects} $V_q$ if
$\Phi_p(W)\cap \im(\Phi_q)\ne \varnothing$.
\enddefinition  

Let us start by looking at two trivial cases: If $p=e$ is a primitive
idempotent, then $V_p$ is a singleton, represented by the
simple module centered at $e$, and if $p=\alpha$ is an arrow $e\rightarrow e'$,
then
$V_p =\Bbb A^n$, where $n$ is the number of arrows $e
\rightarrow e'$ different from $\alpha$. In each of these cases, the
uniserial variety at $p$ consists of a single component.

Observe that, in case the quiver $\Gamma$ of
$\la$ has no double arrows between any two  vertices, the images $\Phi_p(V_p)$
of the uniserial varieties, where
$p$ runs through the paths in $K\Gamma$, yield a disjoint partitioning of the
set of isomorphism classes of uniserial objects in
$\lamod$.  In general, however, the uniserial varieties at $p$ and $q$,
where
$p$ and
$q$ are distinct paths of the same length passing through the same sequence of
vertices, may intersect.  We will see in Section 5 that the uniserial
components at $p$ which intersect $V_q$ and those
components at $q$ which intersect $V_p$ coincide in number and can
be arranged into birationally equivalent pairs. 

Two further problems impose themselves as follow-ups to the preceding
theorem.  One concerns functoriality of the assignment $\la = K\Gamma / I
\mapsto
\bigl(V_p \bigr)_{p\ \text{a path in}\ K\Gamma}$ and, in particular, the
behavior of the family of varieties $V_p$ under algebra isomorphism.
The other is the isomorphism problem for uniserial modules. We defer the
former to [6], and proceed by  tackling the latter.

\head  4.  The isomorphism problem for uniserial modules\endhead 

 Since an obvious necessary
condition for isomorphism of two uniserial modules is a joint mast, we wish to
explicitly describe the equivalence relation on
$V_p$ which partitions $V_p$ into the fibres $\Phi_p^{-1} (U)$, where $U$ runs
through the uniserial modules in the image $\Phi_p(V_p)$.  Theorem
A(III) provides us with a partial answer:  If the path $p$ does not start with
an  oriented cycle, i.e., if $p$ is not of the form $p'c$ where $c$ is an
oriented cycle of positive length, the map
$\Phi_p$ is a bijection, and thus two uniserial modules with mast $p$ are
isomorphic if and only if the corresponding points on the variety $V_p$
coincide; in other words, the points on the variety $V_p$ form a complete
system of isomorphism invariants for the uniserials with mast $p$ in that
case.  In general, injectivity of
$\Phi_p$ may fail, as we know from Example 7.  To fill the resulting gap in
the information on the uniserial modules with mast $p$  stored in $V_p$, we
will construct a system of equations
$S_p(X,Y,Z)$ in the variables
$X_i(\alpha,u)$, $Y_i(\alpha,u)$ (for $(\alpha,u) \detour p$ and $i\in
I(\alpha,u)$) and finitely many variables $Z_i$, which is linear in the $Z_i$
over
$K[X,Y]$ such that, for any pair of points $k,k'\in V_p$, the linear system
$S_p(k,k',Z)$ is consistent if and only if $\Phi_p(k)\cong \Phi_p(k')$. So,
loosely speaking, the family of uniserial modules with mast $p$ can be
identified with the variety
$V_p$ modulo a certain system of linear equations with coefficients in $K$. 
While easy to establish, this system will provide a handy  decision
process for the isomorphism problem on the basis of the varieties
$V_p$ which, in turn, are readily accessible.

To describe the system $S_p(X,Y,Z)$, suppose that the path $p :
e(1)\rightarrow e(l+1)$ has precisely $t$ right subpaths of positive length
ending in the starting vertex $e(1)$ of $p$, say $w_1,\dots,w_t$. Then our
system will have the $t$ linear variables $Z_1,\dots,Z_t$. Start by
considering the following equations ($*$) in $K\Gamma[X,Y,Z]$, one for each
detour $(\alpha,u)$ on $p$:
$$\sum_{i\in I(\alpha, u)}
X_i(\alpha, u) v_i(\alpha, u) \biggl( e(1) + \sum_{j=1}^t Z_j w_j \biggr) = 
 \alpha u \biggl( e(1) + \sum_{j=1}^t Z_j w_j \biggr); \tag{$*$}$$  
here the $v_i(\alpha,u)$ are as in the definition of the substitution
equations for $p$. Now expand both sides of these equations by successively
inserting from the right the substitution equations $\beta v \seq \sum_{i \in
I(\beta, v)} Y_i(\beta, v) v_i(\beta, v)$ for  detours
$(\beta,v)$ on $p$, and the equivalences $q\seq0$ for those paths $q\in
K\Gamma$ which fail to be routes on $p$. As pointed out in  Observation 3 of
Section 3, the equations ($*$) will eventually take on the form
$$\sum_{i \in I(\alpha, u)}
a_i(X,Y,Z) v_i(\alpha, u) = \sum_{i \in I(\alpha, u)} b_i(X,Y,Z) v_i(\alpha,
u)$$   
for suitable polynomials $a_i(X,Y,Z),\ b_i(X,Y,Z)\in K[X,Y,Z]$ which are
uniquely determined by the left-hand and right-hand sides of  equation
($*$). Indeed, this is always the terminal stage of the substitution process,
since each substitution step replaces a path $q$ by a linear combination of
paths of lengths $\ge \len(q)$, all of which have the same endpoint as $q$.
Now collect all of the equations of the form $a_i(X,Y,Z)=
b_i(X,Y,Z)$, $i\in I(\alpha,u)$, arising in this way for arbitrary detours
$(\alpha,u)$ on $p$, and label the resulting system $S_p(X,Y,Z)$. Observe that
this system is polynomial in the $X_j$ and  $Y_j$, linear in the $Z_j$.

\proclaim {Theorem B} For $k,k' \in V_p$, the linear system $S_p(k,k',Z)$ in
$Z= (Z_1,\dots,Z_t)$ is consistent if and only if $\Phi_p(k)\cong
\Phi_p(k')$.\endproclaim

\demo{Proof}  
First suppose that there exists an isomorphism
$f:\Phi_p(k) \rightarrow \Phi_p(k')$ of $\la$-modules, and let
$x$ and $y$ be top elements of $\Phi_p(k)$ and $\Phi_p(k')$,
with the property that $\alpha u x = \sum_{i \in I(\alpha, u)} k_i(\alpha, u)
v_i(\alpha, u) x$ and $\alpha u y = \sum_{i \in I(\alpha, u)} k'_i(\alpha, u)
v_i(\alpha, u) y$, respectively, for any detour $(\alpha,u)$ on $p$ (here we
identify the isomorphism classes $\Phi_p(k)$ and $\Phi_p(k')$ with the
distinguished representatives described in Theorem A).  The equation
$f(x) =e(1)f(x)$ clearly yields
$f(x) =c_0y+ \sum_{j=0}^t c_j w_jy$ with $c_j \in K$ and $c_0 \ne 0$;  without
loss of generality, we may assume $c_0 = 1$.  Let $(\alpha, u)$ be a detour on
$p$ such that $\alpha$ ends in the primitive idempotent $e$, and
consider the  equality inside the $K$-vectorspace
$e \Phi_p(k')$ which results from our isomorphism $f$:  
$$\sum_{i \in I(\alpha, u)}
k_i(\alpha,u) v_i(\alpha, u) \biggl(y + \sum_{j=1}^t c_jw_j y \biggr) =
\alpha u \biggl(y + \sum_{j=1}^t c_j w_jy \biggr). \tag{$\dagger$}$$  Clearly,
the left-hand side of equality ($\dagger$) reduces to the form 
$$\sum_{i\in I(\alpha,u)} a_i(k,k',c) v_i(\alpha,u)y \in e\Phi_p(k'),$$
where $c= (c_1,\dots,c_t)$ and $a_i(X,Y,Z)$ is the
polynomial occurring in the definition of the system $S_p(X,Y,Z)$; this is
clear from the discussion prior to Theorem B. Analogously, the
right-hand side of ($\dagger$) is equal to the element
$$\sum_{i\in I(\alpha,u)} b_i(k,k',c) v_i(\alpha,u)y \in e\Phi_p(k')$$
with
$b_i(X,Y,Z)
\in K[X,Y,Z]$ as above, and consequently equality ($\dagger$) reduces to
$$\sum_{i \in I(\alpha, u)}
a_i(k,k',c) v_i(\alpha, u) y = \sum_{i \in I(\alpha, u)} b_i(k,k',c) v_i(\alpha,
u) y. \tag{$\ddagger$}$$ 
Since the vectors
$v_i(\alpha, u) y$, $i \in I(\alpha, u)$, are $K$-linearly independent,
this shows that the scalars $c_1, \dots, c_t$ satisfy the system $S_p(k,k',Z)$.

Conversely, suppose that the system $S_p(k,k',Z)$ is consistent, and let
$c = (c_1, \dots , c_t) \in K^t$ be a solution.  Moreover, denote the residue
class of the vertex $e(1)$ in $\Phi_p(k)$ by $x$, that in  $\Phi_p(k')$ by
$y$ (again we identify $\Phi_p(k)$ and $\Phi_p(k')$ with the representatives
described in Theorem A).  We wish to show that the assignment
$x\mapsto y +
\sum_{j=1}^t c_j w_j y$ extends to a well-defined $\la$-isomorphism $\Phi_p(k)
\rightarrow \Phi_p(k')$.  It clearly suffices to show that the annihilator
of $x$ in $\la$ is contained in the annihilator of $y + \sum_{j=1}^t c_j
w_j y$, since these two elements are top elements of the uniserials 
$\Phi_p(k)$ and $\Phi_p(k')$, respectively, and the uniserial modules
$\Phi_p(k)$ and $\Phi_p(k')$ have the same length.  Recall that, by the
definiton of
$\Phi_p(k)$, the left annihilator of
$x$ is, as a left ideal of $\la$, generated by the (residue classes in $\la$
of) the elements
$\alpha u -\sum_{i\in I(\alpha, u)} k_i(\alpha, u) v_i(\alpha, u)$, where
$(\alpha,u)$ runs through the detours on $p$, and by all the (residue classes
of) paths $q$ which fail to be routes on
$p$.  That each of the former elements annihilates $y + \sum_{j=1}^t c_j
w_j y$ is guaranteed by the fact that $c= (c_1,\dots,c_t)$ satisfies the system
$S_p(k,k',Z)$. Indeed, from the fact that $c$ satisfies ($\ddagger$), we infer
that $c$ also satisfies equality ($\dagger$) above. So let
$q$ be a non-route on
$p$.  By construction,
$q$ also annihilates $y$; moreover, given any of the subpaths $w_j$ of $p$
ending in $e(1)$,  the composition $q w_j$ is
not a route on $p$ either by Remark (a) following Proposition 2.1.  This
implies that also
$q w_j y = 0$ for all
$j$, and hence that $q(y + \sum_{j=1}^t c_j w_j y) = 0$ as required. \qed
\enddemo

We illustrate our method for solving the isomorphism problem with two
examples.  In the first, $V_p \cong \Bbb A^3$ and each of the fibres of
$\Phi_p$ is a subvariety isomorphic to $\Bbb A^2$, in the second $V_p \cong
\Bbb A^1$ consists of a single fibre.  

\example{Example 9} Let $\la=K\Gamma/I$, where $\Gamma$ is the quiver

\ignore{
$$\xymatrixcolsep{3pc}
\xy\xymatrix{
1 \lloopd{\alpha} \ar[r]<0.5ex>^\beta &2 \ar[l]<0.5ex>^\gamma
}\endxy$$
}

\noindent and $I$ is the ideal in $K\Gamma$ generated by $\alpha^2$,
$\gamma\beta\gamma$, $\gamma\beta\alpha\gamma$. Consider the path $p=
\beta\alpha\gamma\beta\alpha$. The detours on $p$ are $(\beta,e_1)$, $(\beta,
\gamma\beta\alpha)$, and $(\alpha,\alpha)$, yielding the substitution
equations 
$$\beta \seq X_1\beta\alpha +X_2p,\qquad \beta\gamma\beta\alpha \seq X_3p,
\qquad\text{and}\qquad \alpha^2 \seq X_4\gamma\beta\alpha +X_5\alpha\gamma
\beta\alpha.$$
Inserting the substitution equations from the right into a $K$-generating set
for $I^{(6)}$ -- note that 6 is the Loewy length of $\Lambda$ -- gives us
$X_4=X_5=0$, while imposing no conditions on $X_1$, $X_2$, $X_3$. Thus $V_p
\cong \Bbb A^3$. 

To
determine the system $S_p(X,Y,Z)$, we observe that there are precisely three
right subpaths of
$p$ of positive length which end in $e_1$, namely $\alpha$,
$\gamma\beta\alpha$, and
$\alpha\gamma\beta\alpha$; thus $Z = (Z_1, Z_2, Z_3)$. Setting $z = e_1
+Z_1\alpha +Z_2\gamma\beta\alpha +Z_3\alpha\gamma\beta\alpha$ and inserting the
substitution equations
$\beta \seq Y_1\beta\alpha +Y_2p$, $\beta\gamma\beta\alpha \seq Y_3p$ and
$\alpha^2 \seq 0$ repeatedly into the three starting equations (*) 
$\beta z =X_1\beta\alpha z+X_2pz$, $\beta\gamma\beta\alpha z =X_3pz$ and
$\alpha^2z =X_4\gamma\beta\alpha z +X_5\alpha\gamma\beta\alpha z$ as
described ahead of Theorem B, we obtain
$$\align X_1\beta\alpha +Z_2X_1p +X_2p &=Y_1\beta\alpha +Y_2p +Z_1\beta\alpha
+Z_2\beta\gamma\beta\alpha +Z_3p\\
 &=Y_1\beta\alpha +Y_2p +Z_1\beta\alpha
+Z_2Y_3p +Z_3p \endalign$$
from the first of the equations (*), $X_3p =Y_3p$ from the second, and
$0=0$ from the third. Thus the system
$S_p(X,Y,Z)$ is
$$X_1=Y_1+Z_1,\qquad Z_2X_1 +X_2=Y_2+Z_2Y_3+Z_3,\qquad X_3=Y_3$$
in this example. In particular, given any two points $k= (k_1,k_2,k_3)$ and
$k'= (k'_1,k'_2,k'_3)$ in $V_p\cong \Bbb A^3$, the linear system $S_p(k,k',Z)$
in
$Z_1$, $Z_2$, $Z_3$  is consistent if and only if
$k_3=k'_3$. This shows that, up to isomorphism, there is only a one-parameter
family of uniserial modules with mast $p$, the parameter being $k_3\in K$.

Let us translate this information into the graphs of the uniserial modules
with mast $p$. The uniserial modules corresponding to the nonzero values of
$k_3$ all have three different graphs, depending on the choice of top
element.  These graphs are:

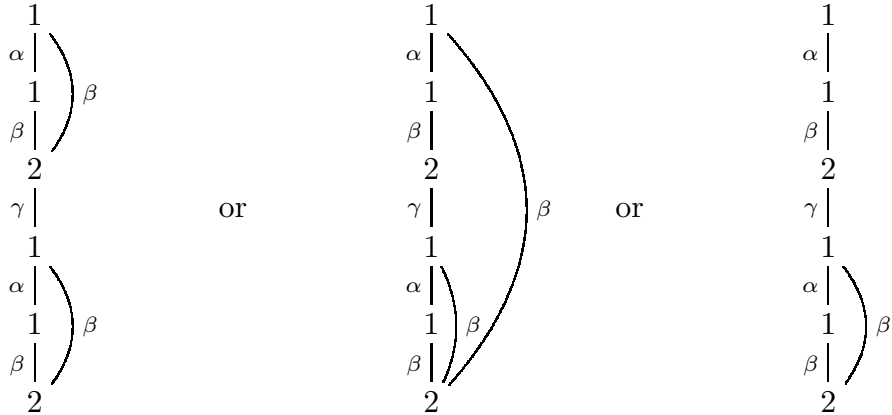
\ignore{
$$\xymatrixrowsep{1pc}
\xy\xymatrix{
1 \edge[d]_\alpha \edge @/^1pc/ [dd]^\beta &&&& 1
\edge[d]_\alpha
\edge @/^2.5pc/ [5,0]^\beta &&&& 1 \edge[d]_\alpha \\
1 \edge[d]_\beta &&&& 1 \edge[d]_\beta &&&& 1 \edge[d]_\beta
\\ 
2 \edge[d]_\gamma && \ar@{}[d]|{\displaystyle{\txt{or}}} && 2
  \edge[d]_\gamma && \ar@{}[d]|{\displaystyle{\txt{or}}} && 2
\edge[d]_\gamma \\  
1 \edge[d]_\alpha \edge @/^1pc/ [dd]^\beta &&&& 1
\edge[d]_\alpha \edge @/^0.65pc/ [dd]^\beta &&&& 1
\edge[d]_\alpha \edge @/^1pc/ [dd]^\beta \\
1 \edge[d]_\beta &&&& 1 \edge[d]_\beta &&&& 1 \edge[d]_\beta
\\ 
2 &&&& 2 &&&& 2 }
\endxy$$
}

\noindent The uniserial module corresponding to the value $k_3=0$ has graphs

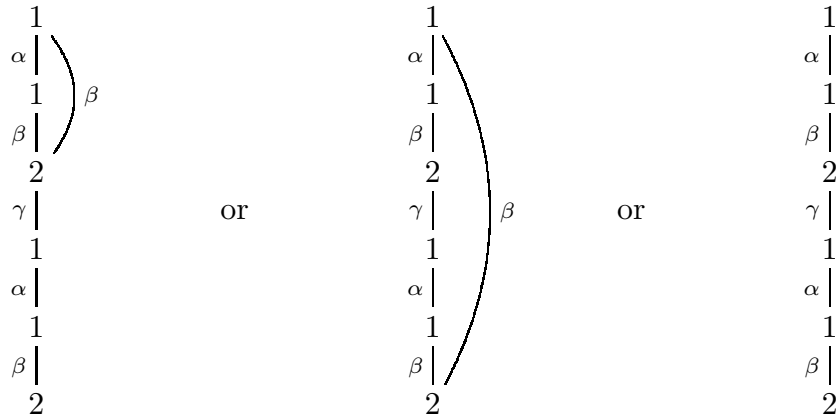
\ignore{
$$\xymatrixrowsep{1pc}
\xy\xymatrix{
1 \edge[d]_\alpha \edge @/^1pc/ [dd]^\beta &&&& 1
\edge[d]_\alpha
\edge @/^1.5pc/ [5,0]^\beta &&&& 1 \edge[d]_\alpha \\
1 \edge[d]_\beta &&&& 1 \edge[d]_\beta &&&& 1 \edge[d]_\beta \\
2 \edge[d]_\gamma && \ar@{}[d]|{\displaystyle{\txt{or}}} && 2
  \edge[d]_\gamma && \ar@{}[d]|{\displaystyle{\txt{or}}} &&
2 \edge[d]_\gamma \\ 
1 \edge[d]_\alpha &&&& 1 \edge[d]_\alpha &&&& 1 \edge[d]_\alpha
\\  
1 \edge[d]_\beta &&&& 1 \edge[d]_\beta &&&& 1 \edge[d]_\beta \\
2 &&&& 2 &&&& 2 }
\endxy$$
}

\noindent again depending on the choice of a top element. \qed\endexample

\example{Example 10} Let $\la=K\Gamma/I$, where $\Gamma$ is the quiver

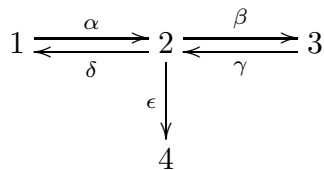
\ignore{
$$\xymatrixcolsep{3pc}
\xy\xymatrix{
1 \ar[r]<0.5ex>^\alpha &2 \ar[r]<0.5ex>^\beta \ar[d]_\epsilon
\ar[l]<0.5ex>^\delta &3 \ar[l]<0.5ex>^\gamma\\
 &4}\endxy$$
}

\noindent and $I$ is the ideal in $K\Gamma$ generated by the relations
$$\epsilon\alpha -\epsilon\gamma\beta\alpha,\quad
\delta\alpha\delta\alpha,\quad
\alpha\delta\alpha\delta,\quad \beta\gamma,\quad \delta\gamma.$$  Observe
that these relations, together with $\epsilon\alpha\delta -
\epsilon\gamma\beta\alpha\delta$ and $\epsilon\alpha\delta\alpha -
\epsilon\gamma\beta\alpha\delta\alpha$ generate $I$ as a left ideal. 
Moreover, consider the mast
$p=
\epsilon\gamma\beta\alpha\delta\alpha$. The detours on
$p$ are $(\beta,\alpha)$, $(\epsilon,\alpha)$, and
$(\epsilon,\alpha\delta\alpha)$. Inserting the corresponding substitution
equations
$$\beta\alpha \seq X_1\beta\alpha\delta\alpha,\qquad \epsilon\alpha \seq
X_2p,\qquad \epsilon\alpha\delta\alpha \seq X_3p$$
into the relations yields $V_p =V(X_2-X_1,\ X_3-1) =\{
(k_1,k_1,1) \mid k_1\in K
\} \cong \Bbb A^1$.

Note
that $w=\delta\alpha$ is the only right subpath of positive length of $p$
which ends in $e_1$, whence the family $Z$ of variables in the system
$S_p(X,Y,Z)$ is reduced to a single one.  Using the above method for
determining the system $S_p(X,Y,Z)$, we obtain, for any two points $k=
(k_1,k_1,1)$ and $k'= (k'_1,k'_1,1)$ in $V_p$:
$$k'_1+Z=k_1,\qquad k'_1+Z=k_1,\qquad 1=1. \tag{$S_p(k,k',Z)$}$$
Since this system is consistent for arbitrary choice of $k,k'\in V_p$, there
is, up to isomorphism, precisely one uniserial left $\la$-module with mast
$p$. The following are all its graphs relative to suitable top elements:

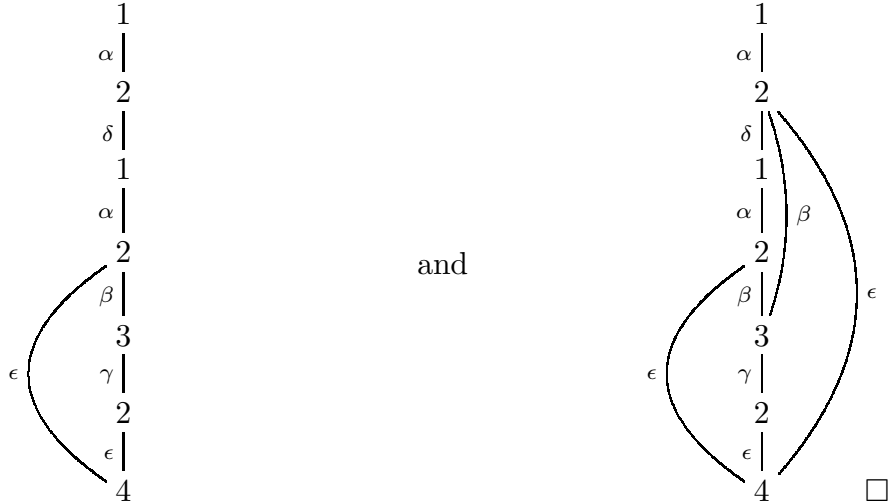
\ignore{
$$\xymatrixrowsep{1pc}
\xy\xymatrix{
1 \edge[d]_\alpha &&&&&& 1 \edge[d]_\alpha \\
2 \edge[d]_\delta &&&&&& 2 \edge[d]_\delta
\edge @/^0.65pc/ [3,0]^\beta \edge @/^2.5pc/ [5,0]^\epsilon \\
1 \edge[d]_\alpha &&&&&& 1 \edge[d]_\alpha \\
2 \edge[d]_\beta \edge @/_2.5pc/ [3,0]_\epsilon &&& \txt{and} &&& 2
\edge[d]_\beta \edge @/_2.5pc/ [3,0]_\epsilon \\
3 \edge[d]_\gamma &&&&&& 3 \edge[d]_\gamma \\
2 \edge[d]_\epsilon &&&&&& 2 \edge[d]_\epsilon \\
4 &&&&&& 4 &\square }
\endxy$$
}
\endexample

We conclude the section with a look at hereditary algebras. As is readily
seen, in that case, all of the varieties $V_p$ are full affine spaces of
`maximum' dimension. We will see that the split hereditary
algebras are actually characterized  by their varieties of uniserial modules,
which answers a question of K\. R\. Fuller.

\proclaim{Proposition C}   The algebra $\la$ is hereditary if and only
if, for each path $p$ in $K\Gamma$, the variety $V_p$ is isomorphic to the
full affine space $\Bbb A^{N(p)}$, where 
$$N(p) = \sum \{|I(\alpha,u)| :
(\alpha,u) \detour p \}$$
and $\Bbb A^0$ stands for a singleton.  

In case these conditions are satisfied, all the maps 
$$\Phi_p: V_p
\longrightarrow \{\text{isomorphism types of uniserials
in\ } \lamod \text{with mast\ } p \}$$ 
are bijections, i.e., $\Phi_p(k) \cong
\Phi_p(k')$ precisely when $k = k'$.  Moreover, for a hereditary algebra
$\la$ over an infinite field $K$, the following statements are equivalent:

{\rm (i)}  There are only finitely many isomorphism types of uniserial
$\la$-modules.

{\rm (ii)}  Given any pair of vertices $e$ and $e'$ in $\Gamma$ and an arrow
$\alpha : e \longrightarrow e'$, the arrow $\alpha$ is the only path from $e$
to $e'$ in $K\Gamma$.

{\rm (iii)}  For any finite sequence of simple $\la$-modules there 
is either no or precisely one uniserial module having this composition
series, depending on whether or not there is a path in $K\Gamma$ which passes
through the corresponding sequence of vertices.\endproclaim

\demo{Proof}  To prove the first equivalence, observe that, for any path
$p$,  the coordinate ring of the variety $V_p$ is $K[X_i(\alpha,u) \mid i \in
I(\alpha,u),\ (\alpha,u) \detour p]$ modulo the ideal generated by all the
polynomials arising from an insertion of the substitution equations into
relations involving routes on $p$.  In particular, this ideal will be nonzero
whenever
$p$ makes a non-trivial appearance in a relation of $\la$, in which case 
$\dim V_p < N(p)$.  Consequently, isomorphism of $V_p$ with $\Bbb A^{N(p)}$ for
all paths $p$ does not allow for any non-trivial relations. Conversely, it is
clear that $V_p\cong \Bbb A^{N(p)}$ for all $p$ when $\la$ is hereditary.

Now suppose that $\la$ is hereditary.  Since the
quiver
$\Gamma$ is acyclic in this case, bijectivity of the maps $\Phi_p$ follows
from part (III) of Theorem A.  In view of the first part of the theorem,
condition (i) is therefore equivalent to `$N(p) = 0$' for all $p$, whence (i)
implies (iii), the converse of this implication being trivial.  But condition
(ii) also translates into the nonexistence of detours on any path $p$ in
$K\Gamma$, that is, into the equality `$N(p) = 0$' for all $p$.  This
completes the proof.      \qed
\enddemo

\head  5. Varieties of uniserials with fixed sequence of composition factors
\endhead
 
A more natural subdivision of the varieties of uniserial modules than that
in terms of masts -- the latter depending a priori on the given
coordinatization of $\la$ -- is in terms of sequences of consecutive
composition factors.  In order to understand the uniserials of
composition length $l+1$ with a fixed sequence $\bigl(S(1), \dots ,
S(l+1)\bigr)$ of simple composition factors, we need to first study the
correlation among the varieties 
$V_p$, where $p$ runs through all paths of length $l$ which pass precisely
through the vertices $e(1), \dots , e(l+1)$, in that order.  In
particular, we need to explore the intersections $\Phi_p(V_p) \cap
\Phi_q(V_q)$, where $p$ and $q$ are two such paths.  In a first easy step, we
will observe that there is a $1-1$ correspondence between the detours on
$p$ and those on $q$ which preserves the cardinalities of the corresponding
index sets $I(\alpha, u)$.  Hence, the varieties $V_p$ and $V_q$ live in the
same affine space $\Bbb A^N$.  Let $D = \Phi_p(V_p) \cap
\Phi_q(V_q)= \im(\Phi_p)\cap \im(\Phi_q)$. It turns
out that the preimages
$\Phi_p^{-1}(D)$ and
$\Phi_q^{-1}(D)$ are Zariski-open in $V_p$ and $V_q$, respectively, and
isomorphic.  (We  do not require varieties to be irreducible and
correspondingly mean by an isomorphism between two varieties a homeomorphism
which has regular coordinate functions in both directions.)  In particular, the
irreducible components of $V_p$ which intersect $V_q$ and those of $V_q$ which
intersect $V_p$ can be paired off into pairs of birationally
equivalent partners. (Recall that, by Definition 8, an irreducible component
$W$ of $V_p$ intersects $V_q$ if $\Phi_p(W)\cap \im(\Phi_q)\ne
\varnothing$.)  The gist of this is that, if we are looking for a set of
representatives of the birational equivalence classes of all the uniserial
varieties of
$\lamod$ at the paths
$p$ running through the above sequence of vertices, we will not lose
information in proceeding as follows:  Let $p_1, \dots , p_t$ be the distinct
paths of length
$l$ passing through the sequence $\bigl(e(1), \dots , e(l+1) \bigr)$.  We start
by determining the irreducible components of $V_{p_1}$.  Then we find the
irreducible components $W$ of $V_{p_2}$ such that $\Phi_{p_2}(W) \cap
\im(\Phi_{p_1}) =
\varnothing$, next the components $W$ of $V_{p_3}$ with the property
that $\Phi_{p_3}(W) \cap \bigl( \im(\Phi_{p_1}) \cup
\im(\Phi_{p_2})
\bigr) =
\varnothing$, and so forth. 
Eventually, this procedure will lead us to a family
of irreducible affine varieties which describes all uniserial $\la$-modules
with a fixed sequence of composition factors and which, up to order and
birational equivalence, is an isomorphism invariant of $\la$. 

\proclaim{Lemma 11}  Suppose that $p$ and $q$ are paths of length $l$ passing
through the same sequence of vertices $\bigl( e(1), \dots , e(l+1) \bigr)$
in the given order.  Then there is a bijection $\rho$ from the set of detours
on
$p$ to the set of detours on $q$ such that $|I \bigl(\alpha,u \bigr)| = |I
\bigl(\rho (\alpha,u) \bigr)|$ for all detours $(\alpha,u) \detour p$.  In
particular, if
$$N = \sum _{(\alpha,u) \detour p} |I(\alpha,u)|,$$ 
then $V_p$ and $V_q$
are both subvarieties of affine $N$-space $\Bbb A^N$ over $K$.
\endproclaim

\demo{Proof} It clearly suffices to focus on the case
where
$p$ and $q$ differ in precisely one arrow; an obvious induction will then
complete the proof.  Say  $p = \alpha_l \cdots \alpha_1$ and $q = \alpha_l
\cdots
\alpha_{r+1} \beta_r \alpha_{r-1} \cdots \alpha_1$.  Let $\rho$ be the
identity on those detours $(\gamma,u)$ on $p$ for which either $\len u < r-1$,
or else $\len u = r-1$ and $\gamma \ne \beta_r$, and set
$\rho(\beta_r,\alpha_{r-1}
\cdots \alpha_1) = (\alpha_r,\alpha_{r-1} \cdots \alpha_1)$ if $\gamma
=\beta_r$.  Any detour of the form 
$(\gamma, \alpha_s \cdots
\alpha_1)$ on $p$ with $s \ge r$, finally, we match up with the detour
$(\gamma, 
\alpha_s \cdots \alpha_{r+1} \beta_r \alpha_{r-1} \cdots \alpha_1)$
on $q$. It is easy to check that this bijection $\rho$ preserves the
cardinalities of the corresponding index sets as claimed. \qed\enddemo

The proof of the following theorem meets with a few technical hurdles which
will, however, cease to play a role in the further development of the subject.

\proclaim{Theorem D} Let $p$ and $q$ be paths of length $l$
passing through the same sequence of vertices, and let $N$ be as in
the preceding lemma.  Moreover, consider the intersection $D =
\Phi_p(V_p)\cap \Phi_q(V_q)$.  Then the preimages $\Phi_p
^{-1}(D)$ and
$\Phi_q ^{-1}(D)$ in $\Bbb A^N$ are isomorphic Zariski open subsets of $V_p$
and $V_q$, respectively.  

More precisely, if $p$ and $q$ differ in exactly $s$ arrows, then there
exist Zariski open subsets $Z_1$ and $Z_2$ of the form
$Z_i = \Bbb A^N\setminus V(X_{i1} \cdots X_{is})$ in $\Bbb A^N$ such that $Z_1
\cap V_p = \Phi_p^{-1}(D)$ and $Z_2 \cap V_q = \Phi_q^{-1}(D)$, together
with an isomorphism $\psi : Z_1 \cap V_p \rightarrow Z_2 \cap V_q$ which
makes the following diagram commutative:

\ignore{
$$\xy\xymatrix{
\Phi_p^{-1}(D) \ar[rr]^\psi \ar[dr]_{\Phi_p} && \Phi_q^{-1}(D)
\ar[dl]^{\Phi_q}\\
 & D
}\endxy$$
}
\endproclaim
 
\demo{Proof} To construct an isomorphism $\psi : \Phi_p^{-1}(D) \rightarrow
\Phi_q^{-1}(D)$, we assume that $D\ne \varnothing$.  We start with a point
$k\in \Phi_p^{-1}(D)$, and set
$U=
\Phi_p(k)$. Moreover, we let $x\in U$ be a top element such that, for each
detour
$(\gamma,u)$ on $p$, we have $\gamma ux =\sum_i k_i(\gamma,u) v_i(\gamma,u)
x$. Then $k$ is the family of
coordinate vectors $(k_i(\gamma,u))_{i\in I(\gamma,u)}$ of the elements
$\gamma ux$ with respect to the basis $(v_i(\gamma,u)x)$, $i\in I(\gamma,u)$,
for the $K$-space $U(\gamma,u) =eJ^{\len(u)+1}U$, where $e$ is the endpoint
of $\gamma$. Let $\rho(\gamma,u) =(\gamma',u')$ where $\rho$ is as in Lemma
11; the lemma then allows us to assume that $I(\gamma',u') =I(\gamma,u)$. As we
will see, the coordinate vector
$\bigl( k'_i(\gamma',u') \bigr)_{i\in I(\gamma',u')}$ of the element
$\gamma'u'x$ relative to the basis $v'_i(\gamma',u')x$, $i\in I(\gamma',u')$,
for
$U(\gamma,u) =U(\gamma',u')$ does not depend on the choice of the element $x$
as above. The assignment $k\mapsto k'= (k'_i(\gamma',u'))_{i\in
I(\gamma',u'), (\gamma',u')\detour q}$ thus yields a well-defined map $\psi :
\Phi_p^{-1}(D) \rightarrow
\Phi_q^{-1}(D)$ by the Corollary to Theorem A. Furthermore, we will find that
the coordinates of $\psi$ are rational functions in the $k_i(\gamma,u)$ which
are defined on the whole domain and which depend only on $p$ and $q$. In
particular, this will show that the coordinate functions of $\psi$ are regular
maps. Once this is established, it will readily follow that $\psi$ is an
isomorphism from
$\Phi_p^{-1}(D)$ to $\Phi_q^{-1}(D)$. Indeed, if
$\psi' : \Phi_q^{-1}(D)
\rightarrow \Phi_p^{-1}(D)$ is the map constructed in complete analogy
to $\psi$, then $\psi'$ has regular coordinate functions by symmetry, and it is
straightforward to check that
$\psi'\psi$ and
$\psi\psi'$ are the identities on
$\Phi_p^{-1}(D)$ and $\Phi_q^{-1}(D)$, respectively.

As in the proof of Lemma 11, we will assume that $p$ and $q$ differ in
precisely one arrow and leave the general case to the
reader. Say
$$p= \alpha_l\cdots \alpha_1\qquad \text{and}\qquad q= \alpha_l\cdots
\alpha_{r+1}\beta_r\alpha_{r-1} \cdots\alpha_1$$
where each $\alpha_i$ is an arrow $e(i)\rightarrow e(i+1)$ and $\beta_r$ is
an arrow $e(r)\rightarrow e(r+1)$. Set $u_0= \alpha_{r-1}\cdots \alpha_1$.
Clearly, we then have $\Phi_p^{-1}(D) = \{k\in V_p \mid k_1(\beta_r,u_0) \ne
0\}$ and $\Phi_q^{-1}(D) = \{k'\in V_q \mid k'_1(\alpha_r,u_0) \ne 0\}$,
where we have chosen our indices so that $v_1(\beta_r,u_0) =\alpha_ru_0$ and
$v'_1(\alpha_r,u_0) =\beta_ru_0$. In other words, if
$Z_1 =
\Bbb A^N\setminus V(X_1(\beta_r,u_0))$ and 
$Z_2 = \Bbb A^N\setminus V(X_1(\alpha_r,u_0))$, then $\Phi_p^{-1}(D) =Z_1\cap
V_p$ and $\Phi_q^{-1}(D) =Z_2\cap V_q$.

To follow the strategy outlined above, we let $(\gamma,u)$ be a detour on
$p$, and treat the following cases separately.

1. $u=u_0$ and $\gamma=\beta_r$. Then $(\gamma',u')
=(\alpha_r,u_0)$.

2. $u=u_0$ and $\gamma\ne\beta_r$. Then $(\gamma',u')= (\gamma,u)$.

3. $\len(u)< \len(u_0)=r-1$. Then $(\gamma',u') =(\gamma,u)$.

4. $\len(u)> \len(u_0)$. If $u=\alpha_s\cdots\alpha_1$, then
$(\gamma',u') =(\gamma, \alpha_s\cdots\alpha_{r+1} \beta_ru_0)$.

Throughout, assume that $I(\gamma,u)$ is a set of natural numbers and that
the paths
$v_i(\gamma,u)$ and
$v'_i(\gamma',u')$ are ordered by length, that is, 
$$\len v_i(\gamma,u) <\len v_{i+1}(\gamma,u)$$
for all $i$, and
similarly for $(\gamma',u')$; then, clearly, $\len v_i(\gamma,u) =\len
v'_i(\gamma',u')$ for all $i\in I(\gamma,u)= I(\gamma',u')$. Again, let $k\in
\Phi_p^{-1}(D) =Z_1\cap V_p$, and identify $\Phi_p(k)$ with a uniserial
module $U$ having top element $x$ such that $\delta vx =\sum k_i(\delta,v)
v_i(\delta,v)$ for each detour $(\delta,v)$ on $p$.

In case 1, we write each path $v'_i(\alpha_r,u_0)$ in the form $u_i\beta_ru_0$,
where $u_i$ is a path of length $\ge 0$; in particular, since
$v'_1(\alpha_r,u_0) =\beta_ru_0$, we have $u_1=e(r+1)$. We infer that
$v_i(\beta_r,u_0) = u_i\alpha_ru_0 =u_iv_1(\beta_r,u_0)$ and compute 
$$\align \alpha_ru_0x &= \sum_{i\in I(\alpha_r,u_0)} k'_i(\alpha_r,u_0)
v'_i(\alpha_r,u_0)x\\
 &= \sum_{i\in I(\alpha_r,u_0)} k'_i(\alpha_r,u_0) u_i \Biggl(\sum_{j\in
I(\beta_r,u_0)} k_j(\beta_r,u_0) v_j(\beta_r,u_0)x \Biggr)\\
 &= \sum_{i\ge 1} k'_i(\alpha_r,u_0) k_1(\beta_r,u_0) v_i(\beta_r,u_0) x
+\sum_{{i\ge 1}\atop {j\ge 2}} k'_i(\alpha_r,u_0) k_j(\beta_r,u_0)
u_iv_j(\beta_r,u_0) x. \endalign$$
Observing that, for $i\ge 1$ and $j\ge 2$, we have
$\len(u_iv_j(\beta_r,u_0))> \len(u_iv_1(\beta_r,u_0))
=\len(v_i(\beta_r,u_0))$ and $\len(u_iv_j(\beta_r,u_0))\ge
\len(v_j(\beta_r,u_0))$, we obtain
$$u_iv_j(\beta_r,u_0)x =\sum_{s\ge \max(i+1,j)} \sigma_{ijs}
v_s(\beta_r,u_0)x$$
for all $i\ge 1$ and $j\ge 2$, where the $\sigma_{ijs}$ are polynomials in
the $k_t(\delta,v)$, for detours $(\delta,v)$ on $p$, which -- aside from
dependence on the indices -- depend only on $p$, $q$ and the detour
considered.

Inserting the second equality into the first yields
$$\alpha_ru_0x =\sum_{i\ge 1} \Biggl( k'_i(\alpha_r,u_0)k_1(\beta_r,u_0)
+\sum_{{1\le s\le i-1}\atop{2\le j\le i}} k'_s(\alpha_r,u_0)k_j(\beta_r,u_0)
\sigma_{sji} \Biggr) v_i(\beta_r,u_0)x.$$
On the other hand, $\alpha_ru_0x =1\cdot v_1(\beta_r,u_0)x$, and since the
elements $v_i(\beta_r,u_0)x$ of $U$ are $K$-linearly
independent, a comparison of coefficients leads us to the system of equations
$k'_1(\alpha_r,u_0)k_1(\beta_r,u_0) =1$ and
$$k'_i(\alpha_r,u_0)k_1(\beta_r,u_0) +\sum_{1\le s\le i-1} k'_s(\alpha_r,u_0)
\Biggl( \sum_{2\le j\le i} k_j(\beta_r,u_0)\sigma_{sji} \Biggr) =0$$
for $i\ge 2$. This system for the `unknowns' $k'_i(\alpha_r,u_0)$ has size
$|I(\beta_r,u_0)| \times |I(\beta_r,u_0)|$ and is  lower triangular, the
scalar $k_1(\beta_r,u_0)$ holding all positions along the main diagonal of
the coefficient matrix; this scalar is nonzero, because we chose $k\in Z_1$.
Consequently, the system is uniquely solvable for the $k'_i(\alpha_r,u_0)$,
and the solution is of the form $k'_i(\alpha_r,u_0)
=\tau_i/(k_1(\beta_r,u_0))^i$, where the
$\tau_i$ are polynomials in the $k_t(\delta,v)$ depending only on $i$, $p$,
$q$ and on the detour $(\alpha_r,u_0)$ on $p$.

We will leave the details of cases 2,3 to the reader, but will carry out case
4. In that case, $\len v'_i(\gamma',u') =\len v_i(\gamma,u) >\len(u')
=\len(u) >\len(u_0)$, and hence $v'_i(\gamma',u') =u_i\beta_ru_0$ for some path
$u_i$ of length $\ge 1$, while $v_i(\gamma,u) =u_i\alpha_ru_0$ for all $i\in
I(\gamma',u') =I(\gamma,u)$.

Again we compute $\gamma'u'x$ in two ways. On one hand,
$$\align \gamma'u'x &= \sum_{i\in I(\gamma',u')} k'_i(\gamma',u')
u_i\beta_ru_0x\\
 &= \sum_{i\in I(\gamma,u)} k'_i(\gamma',u')u_i \Biggl( \sum_{j\in
I(\beta_r,u_0)} k_j(\beta_r,u_0) v_j(\beta_r,u_0)x \Biggr) \\
 &= \sum_{i\in I(\gamma,u)} k'_i(\gamma',u') k_1(\beta_r,u_0) v_i(\gamma,u)x
+\sum_{{i\in I(\gamma,u)}\atop{j\in I(\beta_r,u_0), j\ge 2}} k'_i(\gamma',u')
k_j(\beta_r,u_0) u_iv_j(\beta_r,u_0)x. \endalign$$
Here we use the facts that $I(\gamma',u') =I(\gamma,u)$ and that
$u_iv_1(\beta_r,u_0) =u_i\alpha_ru_0 =v_i(\gamma,u)$. To evaluate the terms
$u_iv_j(\beta_r,u_0)x$ for
$j\ge 2$, we observe that
$$\len(u_iv_j(\beta_r,u_0)) >\len(u_iv_1(\beta_r,u_0))
=\len(u_i\alpha_ru_0) =\len(v_i(\gamma,u))$$
for all $j\ge 2$; moreover, $\len(u_iv_j(\beta_r,u_0)) \ge
\len(v_j(\gamma,u))$, because $u_iv_1(\beta_r,u_0)$, $u_iv_2(\beta_r,u_0)$,
\dots, $u_iv_j(\beta_r,u_0)$ is a sequence of right subpaths of $p$ which
have strictly increasing lengths exceeding $\len(u)$ and which end in the
same vertex as $\gamma$.

Thus, for $j\ge 2$, we have $u_iv_j(\beta_r,u_0)x =\sum_{s\ge \max(i+1,j)}
\sigma_{ijs} v_s(\gamma,u)x$, where the $\sigma_{ijs}$ are again polynomials in
the
$k_t(\delta,v)$ depending solely on the indices, the paths
$p$, $q$, and the detour $(\gamma,u)$. Inserting this information into the
above equality yields
\TagsAsText
$$\multline \gamma'u'x =\\
\sum_{i\in I(\gamma,u)} \Biggl( k'_i(\gamma',u')
k_1(\beta_r,u_0) +\sum_{{1\le s\le i-1}\atop{s\in I(\gamma,u)}}
k'_s(\gamma',u')
\sum_{{j\in I(\beta_r,u_0)}\atop{2\le j\le i}} k_j(\beta_r,u_0) \sigma_{sji}
\Biggr) v_i(\gamma,u)x. \endmultline \tag{I}$$

On the other hand, we can write $u=w\alpha_ru_0$ for a suitable path $w$ of
length $\ge 0$. Then $u' =w\beta_ru_0$, and we obtain
$$\gamma'u'x =\gamma w\beta_ru_0x =\sum_{j\in I(\beta_r,u_0)} k_j(\beta_r,u_0)
\gamma wv_j(\beta_r,u_0)x.$$
Now each of the paths $\gamma wv_j(\beta_r,u_0)$ is longer than $u$ and
has the same endpoint as $\gamma$. Therefore
$$\gamma wv_j(\beta_r,u_0)x =\sum_{i\in I(\gamma,u)} \tau_{ji}
v_i(\gamma,u)x,$$
where the $v_i(\gamma,u)$ are suitable polynomials in the $k_t(\delta,v)$. This
gives us
$$\gamma'u'x =\sum_{i\in I(\gamma,u)} \Biggl( \sum_{j\in I(\beta_r,u_0)}
k_j(\beta_r,u_0) \tau_{ji} \Biggr) v_i(\gamma,u)x, \tag{II}$$
and a comparison of coefficients of the $K$-linearly independent elements
$v_i(\gamma,u)x
\in U$ in equations (I), (II) once more yields a square system of equations for
the
$k'_i(\gamma',u')$. Again the coefficient matrix is lower triangular, and all
diagonal positions are occupied by $k_1(\beta_r,u_0)$. Since $k_1(\beta_r,u_0)
\ne 0$ by the choice of $k$, the system has a unique solution, expressing the
$k'_i(\gamma',u')$ as polynomials in the $k_t(\delta,v)$ divided by powers of
$k_1(\beta_r,u_0)$.
\qed\enddemo

Let $W_1,\dots,W_m$ be a full set of representatives for the birational
equivalence classes of {\it all} the irreducible components occurring in
the varieties $V_p$, where $p$ runs through the sequence of vertices
$(e(1),\dots,e(l+1))$. Since, in particular, $W_i$ is not birationally
equivalent to $W_j$ for $i\ne j$,  Theorem C tells us that, for any choice of
$i\ne j$ and
$p\ne q$ with $W_i\subseteq V_p$ and $W_j\subseteq V_q$, the intersection
$\Phi_p(W_i)\cap \Phi_q(W_j)$ is empty;
in other words, among the irreducible varieties $W_i$ selected above, any two
 corresponding to different masts are disjoint. As announced in the beginning
of the section, we can therefore determine the set $\{W_1,\dots,W_m\}$ as
follows: Let
$p_1,\dots,p_t$ be the distinct paths of length $l$ passing through the
sequence of vertices
$(e(1),\dots,e(l+1))$. Start by listing the irreducible components $W\subseteq
V_{p_1}$; then find the irreducible components $W\subseteq V_{p_2}$ with
$\Phi_{p_2}(W)\cap \Phi_{p_1}(V_{p_1})
=\varnothing$; next the components $W\subseteq V_{p_3}$ with
$\Phi_{p_3}(W)\cap (\Phi_{p_1}(V_{p_1})\cup \Phi_{p_2}(V_{p_2})) =\varnothing$;
and so forth. In particular, Theorem C  guarantees that the result does not
depend on the ordering of the paths $p_1,\dots,p_t$. In other words, the set
of birational equivalence classes obtained will be invariant under
permutation of the $p_i$.

We interrupt the theory to give two examples. In each case, $p$ and $q$ are
two paths of the same length passing through the same sequence of vertices.
In the first example, $V_p$ and $V_q$ are irreducible with $\Phi_p(V_p)\cap
\Phi_q(V_q) \ne\varnothing$, but $V_p\not\cong V_q$; Theorem D guarantees that
$V_p$ and $V_q$ are birationally equivalent in that case. In the second, again
$\im(V_p)\cap \im(V_q)\ne \varnothing$, but this time
$V_p$ is irreducible while $V_q$ is not; in this situation, Theorem D
guarantees that $V_p$ is birationally equivalent to an irreducible component
of $V_q$.

\example{Example 12} Let $\Gamma$ be the quiver

\ignore{
$$\xymatrixrowsep{1pc}\xymatrixcolsep{3pc}
\xy\xymatrix{
1 \ar[r]<0.5ex>^{\alpha'} \ar[r]<-0.5ex>_\alpha &2
\ar[r]<0.5ex>^{\beta'} \ar[r]<-0.5ex>_\beta &3
}\endxy$$
}

\noindent and $\Lambda =K\Gamma/I$, where $I$ is generated by $\beta'\alpha'
-\beta\alpha$. If $p$ is equal to $\beta\alpha$ or $\beta'\alpha'$ and
$q$ to $\beta'\alpha$ or $\beta\alpha'$, then
$V_p =V(X_1X_2-1)$ is isomorphic to the punctured line, while $V_q=
V(X_1-X_2)$ is isomorphic to the full affine line over $K$. Thus $V_p\not\cong
V_q$. Note, however, that $\Phi_p(V_p)\cap
\Phi_q(V_q) \ne\varnothing$. \qed
\endexample

\example{Example 13} Let $\Gamma$ be the quiver

\ignore{
$$\xymatrixrowsep{1pc}\xymatrixcolsep{3pc}
\xy\xymatrix{
1 \ar[r]<0.5ex>^{\alpha'} \ar[r]<-0.5ex>_\alpha &2
\ar[r]<0.5ex>^{\beta'} \ar[r]<-0.5ex>_\beta &3
\ar[r]<0.5ex>^{\gamma'} \ar[r]<-0.5ex>_\gamma &4
\ar[r]<0.5ex>^{\delta'} \ar[r]<-0.5ex>_\delta &5
\ar[r]<0.5ex>^{\epsilon'} \ar[r]<-0.5ex>_\epsilon &6
}\endxy$$
}

\noindent and $\Lambda =K\Gamma/I$, where $I$ is the ideal generated by the
relations
$$\epsilon'\delta'\gamma\beta\alpha
-\epsilon\delta\gamma'\beta'\alpha',\qquad
\epsilon\delta\gamma\beta'\alpha -\epsilon\delta\gamma\beta\alpha',\qquad
\epsilon\delta\gamma'\beta\alpha -\epsilon\delta\gamma\beta\alpha',\qquad
\epsilon\delta'\gamma\beta\alpha -\epsilon'\delta\gamma\beta\alpha.$$
If $p=\epsilon\delta\gamma\beta\alpha$ and $q=
\epsilon'\delta\gamma\beta'\alpha$, then 
$$V_p =V(X_4X_5-X_1X_2X_3,\ X_2-X_1,\ X_3-X_1,\ X_4-X_5) \cong V(Y^2-X^3)$$
is the standard cusp:

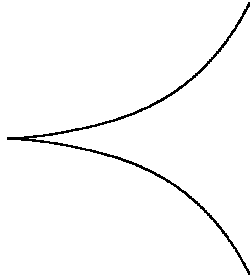
\ignore{
$$\xy (30,15);(30,-15)
\curve{(25,5)&(10,0)&(0,0)&(10,0)&(25,-5)}
\endxy$$
}

\noindent In particular, $V_p$ is irreducible. On the other hand,
$$\align V_q &=V(X_2,X_5)\cup V(X_4X_2-X_5X_3X_1,\ X_1X_2-1,\ X_3-X_1,\
X_5X_4-1)\\
 &\cong \Bbb A^3 \cup V(X_4^2-X_1^3,\ X_1X_2-1,\ X_5X_4-1), \endalign$$ 
the latter component being  isomorphic to the cusp from which the singular
point has been deleted, or, in other words, isomorphic to the punctured line.
Note that, in this case, still $\Phi_p(V_p)\cap
\Phi_q(V_q) \ne\varnothing$, but
$\Phi_p(V_p)\cap \Phi_q(V(X_2,X_5)) =\varnothing$.
\qed
\endexample

As repeatedly announced, the birational equivalence classes of the irreducible
components of the varieties
$V_p$, where $p$ runs through the paths of length $l$  through a fixed sequence
$(e(1),\dots,e(l+1))$ of vertices, is an isomorphism invariant of the algebra
$\Lambda$. To see this, we will require a technical theorem which describes the
behavior of the varieties $V_p$ under change of coordinatization. The proof of
this theorem constitutes a rather lengthy detour from our main train of
thought; in particular, it is based on an extension of our conceptual
framework. Therefore, it will be carried out separately in [6].

\proclaim{Theorem E} Given are two coordinatizations of $\la$, namely $\Lambda
=K\Gamma/I
\cong K{\widehat \Gamma}/{\widehat I}$, where the quiver $\Gamma$ is based on
the primitive idempotents $e_1,\dots,e_n$ as before, and $\widehat \Gamma$
(necessarily isomorphic to $\Gamma$ as a directed graph) is based on
primitive idempotents ${\widehat e}_1,\dots,{\widehat e}_n$ of $\Lambda$,
ordered in such a way that each $e_i$ is congruent, modulo $J$, to the image of
$\widehat e_i$ under our isomorphism.

If $\widehat p$ is a path of length $l$ in $K\widehat \Gamma$ passing through
the sequence of vertices $(\widehat e(1),\dots,\allowmathbreak
\widehat e(l+1))$, then each irreducible component of $V_{\widehat p}$ is
birationally equivalent to a component of some variety $V_p$, where $p$ is a
path of length $l$ in
$K\Gamma$ passing through the sequence $(e(1),\dots,e(l+1))$. Here
$V_{\widehat p}$ is the variety of $\widehat p$ relative to the coordinates
$\widehat \Gamma$ and $\widehat I$, and $V_p$ the variety of $p$ relative to
$\Gamma$ and $I$.
\qed\endproclaim

\proclaim{Theorem F and Definition} Fix a sequence $\Bbb S =(S(1),\dots,
S(l+1))$ of simple left $\Lambda$-modules and let $(e(1),\dots,e(l+1))$ be
the corresponding canonical sequence of primitive idempotents in a fixed
coordinatization $\Lambda =K\Gamma/I$. Then the set of birational equivalence
classes of the irreducible components of the non-empty varieties $V_p$, where
$p$ runs through the paths of length $l$ in $K\Gamma$ passing through
$(e(1),\dots,e(l+1))$ in that order, is independent of the coordinatization.
In other words, the set $V_{\Bbb S}$ consisting of the irreducible components
of the  non-empty varieties
$V_p$ as above is uniquely determined by the isomorphism type of $\la$, up to
birational equivalence. We call $V_{\Bbb S}$ the set of
\underbar{irreducible uniserial varieties of
$\Lambda$ at
$\Bbb S$}.

Moreover, if $V_{\Bbb S} =\{W_1,\dots,W_m\}$, we call the number
$\max_{i\le m} \dim W_i$ the
\underbar{uniserial} \underbar{dimension of $\Lambda$ at $\Bbb S$} and denote
it by
$\dim V_{\Bbb S}$. The \underbar{uniserial genus of $\la$ at $\Bbb S$} is
defined to be the maximum of the numbers 
$\genus(W_1),\dots, \genus(W_m)$, and is denoted by $\genus V_{\Bbb S}$. By the
preceding paragraph, both $\dim V_{\Bbb S}$ and $\genus V_{\Bbb S}$ are
uniquely determined by the isomorphism class of $\la$.
\qed\endproclaim

In Example 12 above, set $\Bbb S= (S_1,S_2,S_3)$. Then $V_{\Bbb S} =\{\Bbb
A^1\}$,  $\dim V_{\Bbb S} =1$, and $\genus V_{\Bbb S} = 0$. In Example 13,
let $\Bbb S= (S_1,\dots,S_6)$. One can readily supplement the computations
given to obtain
$V_{\Bbb S} =\{C, \Bbb A^3\}$, where $C$ is the cusp, and therefore $\dim
V_{\Bbb S} =3$ while, again,
$\genus V_{\Bbb S} =0$. An example of non-trivial uniserial genus can be
found in Section 6 (Example 14).

\head 6. Realization of arbitrary varieties as varieties $V_{\Bbb S}$\endhead

The aim of this section is to show that each affine variety $V$ is isomorphic
to a variety of uniserial modules over a suitable finite dimensional
algebra
$\la =K\Gamma/I$, under the additional
restrictions that $\Gamma$ be acyclic and without double arrows. More
precisely, we will realize our given variety $V$ as a variety $V_p$ for some
path $p$ in $\Gamma$. Observe that, due to the absence of double arrows,
$V_p$ can be identified with
$V_{\Bbb S}$ where $\Bbb
S$ is the sequence of simple modules corresponding to the consecutive
vertices on $p$. Thus the set of irreducible components of an arbitrary
variety
$V$ can be realized as a $V_{\Bbb S}$.

By Theorem
A(III), the first of the above restrictions on
$\Gamma$ guarantees that
$\Phi_p$ is bijective, meaning that there is a 1--1 correspondence between
the points of
$V_p$ and the uniserial left $\la$-modules with mast $p$.
The second condition on $\Gamma$ ensures that all of the varieties $V_q$
where $q$ runs through the paths in $K\Gamma$, are actually isomorphism
invariants of
$\la$ [6]. In particular, our theorem thus
provides a taste of the opulence and complexity of uniserial representations
of finite dimensional algebras.

\proclaim{Theorem G} Given a field $K$ and any affine algebraic variety $V$
over $K$, there exists a finite dimensional path algebra modulo relations,
$\la= K\Gamma/I$, together with a path $p$ in $K\Gamma$, such
that
$V\cong V_p$. We can, moreover, choose the quiver
$\Gamma$ to be acyclic and without double arrows.\endproclaim

\demo{Proof} Say $V= V(f_1,\dots,f_M)$, where the $f_i$ are polynomials in
$K[X_1,\dots,X_m]$ for some $m\ge 1$. In a preliminary step, we construct an
affine variety
$V'$ isomorphic to $V$ such that $V'$ is the vanishing set of polynomials
$f'$ in a certain polynomial ring $K[X_{ij}]$ which have the property that none
of the variables
$X_{ij}$ occurs in a power higher than 1 in any monomial.

To that end, let $d_i= \max_{1\le j\le M} \deg_{X_i}(f_j)$ for $1\le i\le m$,
and introduce new variables $X_{11}$, $X_{12}$, \dots, $X_{1d_1}$, $X_{21}$,
$X_{22}$, \dots, $X_{2d_2}$, \dots, $X_{m1}$, \dots, $X_{md_m}$. For each
$s\in \{1,\dots,M\}$, construct a polynomial $f'_s\in K[X_{ij} \mid 1\le i\le
m,\ 1\le j\le d_i]$ as follows: Replace any monomial $X_1^{r_1}\cdots
X_m^{r_m}$ occurring nontrivially in $f_s$ by
$$(X_{11}\cdots X_{1r_1})(X_{21}\cdots X_{2r_2})\cdots (X_{m1}\cdots
X_{mr_m})$$
in $K[X_{ij}]$; this is possible, since $r_i\le d_i$ for $1\le i\le m$ by
construction. Clearly the total degree of this new monomial in the $X_{ij}$ is
the same as that of $X_1^{r_1}\cdots
X_m^{r_m}$. Define
$$V'= V(f'_1,\dots, f'_M,\ X_{is}-X_{it} \mid 1\le i\le m,\ 1\le s,t\le
d_i).$$
That $V\cong V'$ is then an obvious consequence of our construction, and
hence we may assume, without loss of generality, that $V=V'$; in other words,
we may assume that
$d_i\le 1$ for all $i\in \{1,\dots,m\}$. This concludes the preliminary step.

We now let $\Gamma$ be the quiver

\ignore{
$$\xymatrixcolsep{1.75pc}
\xy\xymatrix{
1 \ar[r]_{\alpha_1} \ar@/^1pc/[rr]^{\gamma_1} &2
\ar[r]_{\beta_1} &3
\ar[r]_{\alpha_2} \ar@/^1pc/[rr]^{\gamma_2} &4 
\ar[r]_{\beta_2} &5
\ar[r]_-{\alpha_3} &\ar@{}[r]|{\displaystyle{\cdots}}
&\ar[r]_-{\beta_{m-1}} &2m-1
\ar[r]_-{\alpha_m} \ar@/^1pc/[rr]^{\gamma_m} &2m 
\ar[r]_-{\beta_m} &2m+1
}\endxy$$
}

\noindent and set $p=\beta_m\alpha_m \cdots \beta_2\alpha_2\beta_1\alpha_1$.
To define suitable relations in $K\Gamma$, we write each $f_j$ in the form
$f_j=
\sum_{A\in \Cal P} c_j(A) \prod_{i\in A} X_i$, where $\Cal P$ is the
power set of $\{1,\dots,m\}$ and $c_j(A)\in K$; this is possible by the
preceding paragraph. For each $A\in \Cal P$, define a path $p(A)$ in
$K\Gamma$ as $p(A) =p_m(A)\cdots p_1(A)$, where the $p_i(A)$ are defined by
$$p_i(A) =\cases \gamma_i\quad &\text{if}\ i\in A\\
\beta_i\alpha_i\quad &\text{if}\ i\notin A\endcases .$$
Note that each $p(A)$ is a path from 1 to $2m+1$. Let
$r_j= \sum_{A\in \Cal P} c_j(A)p(A)$ for $1\le j\le M$, and observe
that each $r_j$ is a relation in $K\Gamma$. Finally,  define $\la
=K\Gamma/I$, where
$I\subseteq K\Gamma$ is the ideal generated by $r_1,\dots,r_M$, and note that
$I^{(L)} =I^{(2m+1)} =I$ is generated by the $r_j$ as a $K$-space. We now
verify that
$V\cong V_p$ as follows.

Clearly, there are no paths in $K\Gamma$ starting in 1 which fail to be
routes on $p$, and the detours on $p$ are precisely the pairs $(\gamma_i,
\beta_{i-1}\alpha_{i-1}\cdots \beta_1\alpha_1)$ for $1\le i\le m$; here
$\beta_0\alpha_0$ stands for the primitive idempotent $e_1$ identified with
the vertex 1. Moreover, for each detour $(\gamma,u)$ on $p$, there exists
precisely one right subpath $v(\gamma,u) =v_1(\gamma,u)$ of $p$ longer than
$u$ and ending in the same vertex as $\gamma$, namely $v(\gamma_i,
\beta_{i-1}\alpha_{i-1}\cdots \beta_1\alpha_1) =\beta_i\alpha_i \cdots
\beta_1\alpha_1$. Writing $X_i$ for the variable $X_1(\gamma_i,
\beta_{i-1}\alpha_{i-1}\cdots \beta_1\alpha_1)$, we thus obtain the
substitution equations
$$\align \gamma_1e_1 &\seq X_1\beta_1\alpha_1\\
\gamma_2\beta_1\alpha_1 &\seq X_2\beta_2\alpha_2\beta_1\alpha_1\\
 &\  \vdots\\
\gamma_m\beta_{m-1}\alpha_{m-1} \cdots\beta_1\alpha_1
&\seq X_m\beta_m\alpha_m \cdots\beta_1\alpha_1 =X_mp\endalign$$
Inserting these successively into the relations $r_j$ from the right yields,
after at most $d$ steps, the equations $r_j\seq f_j(X_1,\dots,X_m)p$, this
being an immediate consequence of our construction. Therefore $V_p
=V(f_1,\dots,f_M) =V$ as desired.
\qed\enddemo

Note that the proof of Theorem G is constructive. We conclude by applying the
pertinent algorithm to obtain an algebra with a uniserial variety of positive
genus.

\example{Example 14}  Let $V= V(X_2^2- X_1(X_1^2-1))$ be the elliptic
curve with $\Bbb R$-graph

\ignore{
$$\xy (0,0);(0,0)
\curve{(0,2)&(-2,7)&(-10,7)&(-10,-7)&(-2,-7)&(0,-2)}
\endxy
\hskip1.0truecm\xy (35,15);(35,-15)
\curve{(30,10)&(18,7)&(10,4)&(10,-4)&(18,-7)&(30,-10)}
\endxy$$
}

\noindent Introduce the new variables $X_{11}, X_{12}, X_{13}, X_{21},
X_{22}$ and note that the variety
$$V' =V(X_{21}X_{22} -X_{11}X_{12}X_{13}+X_{11},\ X_{11}-X_{12},\
X_{11}-X_{13},\ X_{21}-X_{22})$$
is isomorphic to $V$. The following change of indexing makes it easier to
follow the pattern described in the proof of Theorem F:
$$V'= V(X_5X_4-X_3X_2X_1+X_1,\ X_1-X_2,\ X_1-X_3,\ X_4-X_5).$$
Accordingly, we consider the
quiver

\ignore{
$$\xy\xymatrix{
1 \ar[r]_{\alpha_1} \ar@/^1pc/[rr]^{\gamma_1} &2
\ar[r]_{\beta_1} &3
\ar[r]_{\alpha_2} \ar@/^1pc/[rr]^{\gamma_2} &4 
\ar[r]_{\beta_2} &5
\ar[r]_{\alpha_3} \ar@/^1pc/[rr]^{\gamma_3} &6
\ar[r]_{\beta_3} &7 
\ar[r]_{\alpha_4} \ar@/^1pc/[rr]^{\gamma_4} &8
\ar[r]_{\beta_4} &9
\ar[r]_{\alpha_5} \ar@/^1pc/[rr]^{\gamma_5} &10 
\ar[r]_{\beta_5} &11
}\endxy$$
}

\noindent Let $p= q_5q_4q_3q_2q_1$, where $q_i=\beta_i\alpha_i$ for $1\le i\le
5$, and define $\la =K\Gamma/I$, where $I$ is the ideal generated by the
following relations:
$$\gather \gamma_5\gamma_4q_3q_2q_1 -q_5q_4\gamma_3\gamma_2\gamma_1
+q_5q_4q_3q_2\gamma_1,\hskip0.4truein q_5q_4q_3q_2\gamma_1
-q_5q_4q_3\gamma_2q_1,\\
 q_5q_4q_3q_2\gamma_1 -q_5q_4\gamma_3q_2q_1,\hskip0.4truein
q_5\gamma_4q_3q_2q_1 -\gamma_5q_4q_3q_2q_1.\endgather$$ 
Then $V_p\cong V'$, as
substantiated in the proof of the theorem. In particular, if $\Bbb S=
(S(1),\dots,S(11))$, we obtain $\dim V_{\Bbb S} =\genus V_{\Bbb S} =1$.
\qed\endexample

\remark{Remark} Note that, given an affine variety $V$, the $K$-dimension of
the algebra $\la$ over which $V$ is realized as a variety of uniserial
modules increases steeply as the number of variables and their exponents in a
description of $V$ grow. Indeed, if $\Gamma(m)$ is the quiver

\ignore{
$$\xymatrixcolsep{1.75pc}
\xy\xymatrix{
1 \ar[r]_{\alpha_1} \ar@/^1pc/[rr]^{\gamma_1} &2
\ar[r]_{\beta_1} &3
\ar[r]_{\alpha_2} \ar@/^1pc/[rr]^{\gamma_2} &4 
\ar[r]_{\beta_2} &5
\ar[r]_-{\alpha_3} &\ar@{}[r]|{\displaystyle{\cdots}}
&\ar[r]_-{\beta_{m-1}} &2m-1
\ar[r]_-{\alpha_m} \ar@/^1pc/[rr]^{\gamma_m} &2m 
\ar[r]_-{\beta_m} &2m+1
}\endxy$$
}

\noindent then $\dim_K K\Gamma(1) =7$, and $\dim_K K\Gamma(m+1) =4\dim_K
K\Gamma(m)+3$. If one renounces the requirement that the quiver $\Gamma$ of
the algebra realizing $V$ be without double arrows, one can significantly
curb the growth of this dimension in terms of the number of variables
involved. Indeed, the role played by the quiver $\Gamma(m)$ in the proof of
Theorem G can  be taken over by the quiver

\ignore{
$$\xymatrixcolsep{3pc}
\xy\xymatrix{
1 \ar[r]<0.5ex>^{\gamma_1} \ar[r]<-0.5ex>_{\alpha_1}
&2 \ar[r]<0.5ex>^{\gamma_2} \ar[r]<-0.5ex>_{\alpha_2}
&3 \ar[r]<0.5ex>^-{\gamma_3} \ar[r]<-0.5ex>_-{\alpha_3}
&\ar@{}[r]|{\displaystyle{\cdots}}  & \ar[r]<0.5ex>^-{\gamma_{m-1}}
\ar[r]<-0.5ex>_-{\alpha_{m-1}} &m \ar[r]<0.5ex>^-{\gamma_m}
\ar[r]<-0.5ex>_-{\alpha_m} &m+1 }\endxy$$
}

\noindent if one allows for double arrows.
\endremark


\Refs

\ref\no 1\by M. Auslander and I. Reiten \paper Applications of
contravariantly finite subcategories \jour Advances in Math. \vol 86 \yr 1991
\pages 111-152\endref

\ref\no 2\by M. Auslander, I. Reiten, and S. O. Smal\o \book Representation
Theory of Artin Algebras \publaddr Cambridge--New York--Melbourne \yr 1995
\publ Cambridge Univ. Press\endref 

\ref\no 3\by M. Auslander and S. O. Smal\o \paper Preprojective modules over
artin algebras \jour J. Algebra \vol 66 \yr 1980 \pages 61-122\endref

\ref\no 4\by W. D. Burgess and B. Zimmermann Huisgen \paper Approximating
modules by modules of finite projective dimension \jour J. Algebra
\toappear \endref

\ref\no 5\by R. Hartshorne \book Algebraic Geometry \publaddr Berlin \yr 1977
\publ Springer-Verlag\endref

\ref\no 6\by B. Zimmermann Huisgen \paper The geometry of uniserial
representations of finite dimensional algebras II \finalinfo in
preparation\endref

\ref\no 7\by B. Zimmermann Huisgen \paper The geometry of uniserial
representations of finite dimensional algebras III \finalinfo preprint\endref

\ref\no 8\by C. U. Jensen and H. Lenzing \book Model Theoretic Algebra
\bookinfo Algebra, Logic and Applications Series 2 \publ Gordon and Breach
\yr 1989 \publaddr Amsterdam\endref

\ref\no 9\by D. Mumford \book Algebraic Geometry I, Complex Projective
Varieties \publaddr Berlin \yr 1976 \publ Springer-Verlag\endref

\endRefs
\enddocument